\documentclass[11pt,reqno,twoside]{amsart}
\usepackage{amsfonts,amsmath,amssymb}
\usepackage{mathrsfs,mathtools}
\usepackage{enumerate}
\usepackage{hyperref}
\usepackage{esint}
\usepackage{graphicx}
\usepackage{bm}
\usepackage{commath}
\usepackage{esint}
\usepackage{placeins}
\usepackage{flafter}
\usepackage{tikz}

\usepackage[latin1]{inputenc}

\DeclareMathAlphabet{\mathpzc}{OT1}{pzc}{m}{it}
\numberwithin{equation}{section}
\makeatletter
\def\eqnarray{\stepcounter{equation}\let\@currentlabel=\theequation
\global\@eqnswtrue
\tabskip\@centering\let\\=\@eqncr
$$\halign to \displaywidth\bgroup\hfil\global\@eqcnt\z@
  $\displaystyle\tabskip\z@{##}$&\global\@eqcnt\@ne
  \hfil$\displaystyle{{}##{}}$\hfil
  &\global\@eqcnt\tw@ $\displaystyle{##}$\hfil
  \tabskip\@centering&\llap{##}\tabskip\z@\cr}
\def\endeqnarray{\@@eqncr\egroup
      \global\advance\c@equation\m@ne$$\global\@ignoretrue}

\newtheorem{theorem}{Theorem}[section]

\newtheorem{corollary}[theorem]{Corollary}
\newtheorem{definition}[theorem]{Definition}

\newtheorem{proposition}[theorem]{Proposition}
\newtheorem{assumption}[theorem]{Assumption}

\numberwithin{equation}{section}

\def\Omc{\mathbb{R}^N\setminus\Omega}

\def\RR{{\mathbb{R}}}
\def\NN{{\mathbb{N}}}

\def\Om{\Omega}

\def\ga{\gamma}

\def\pOm{\partial\Omega}

\newcommand{\dv}[1]{\,{\mathrm d}#1}
\def\veps{\varepsilon}
\def\R{\mathbb{R}}
\def\O{\Omega}

\def\bxi{{\boldsymbol{\xi}}}
\def\bth{{\boldsymbol{\theta}}}
\def\bu{{\bf{u}}}

\usepackage[textsize=small]{todonotes}
\usepackage[dvips,bottom=1.4in,right=1in,top=1in, left=1in]{geometry}

\thanks{This work is partially supported by NSF grants DMS-1818772, DMS-1913004 
and the Air Force Office of Scientific Research under Award NO: FA9550-19-1-0036.  The work of MW is partially supported by the AFOSR under Award NO:  FA9550-18-1-0242 and by US Army Research Office (ARO) under Award NO: W911NF-20-1-0115. \\
DISTRIBUTION STATEMENT A. Approved for public release. Distribution is unlimited.}



\begin{document}
\title{Optimal Control, Numerics, and Applications of Fractional PDEs}

\author[H.~Antil]{Harbir Antil}
\address{H. Antil,  T.~Brown,  A.~Onwunta,  D. Verma, and M.~Warma. The Center for Mathematics and Artificial Intelligence
(CMAI) and Department of Mathematical Sciences, George Mason University,
Fairfax, VA 22030, USA \\ 
ORCID: 0000-0002-6641-1449}
\email{hantil@gmu.edu, tbrown62@gmu.edu, aonwunta@gmu.edu, dverma2@gmu.edu, mwarma@gmu.edu}

\author[T.~Brown]{Thomas S. Brown}

\author[R.~Khatri]{Ratna Khatri}
\address{R. Khatri. Optical Sciences Division, U.S. Naval Research Laboratory, Washington, DC 20375, USA \\
ORCID: 0000-0003-0931-4025}
\email{ratna.khatri@nrl.navy.mil}

\author[A.~Onwunta]{Akwum Onwunta}

\author[D.~Verma]{Deepanshu Verma}

\author[M.~Warma]{Mahamadi Warma}

\begin{abstract}
This article provides a brief review of recent developments on two nonlocal operators: fractional Laplacian
and fractional time derivative. We start by accounting for several applications of these operators in 
imaging science, geophysics, harmonic maps and deep (machine) learning. 
Various notions of solutions to linear fractional elliptic equations are provided and
numerical schemes for fractional Laplacian and fractional time derivative are discussed. Special emphasis
is given to exterior optimal control problems with a linear elliptic equation as constraints. In addition, optimal control 
problems with interior control and state constraints are considered. We also provide a discussion on fractional 
deep neural networks, which is shown to be a minimization problem with fractional in time ordinary differential
equation as constraint. The paper concludes with a discussion on several open problems. 
\end{abstract}

\maketitle
\tableofcontents

\section{Introduction and applications of fractional operators}

Fractional (nonlocal) models have received a significant amount of attention during the recent years.
This can be attributed to their ability to account for long range interactions and their flexibility 
in being applicable to non-smooth functions. Motivated by these facts and several applications,
this research area has attracted a significant amount of attention on 
both theoretical and computational fronts. The goal of this article is to discuss some of these 
developments. This article is meant to provide an overview, largely motivated by authors' own 
research, and it is not meant to be completely exhaustive. 
We begin by discussing several applications of fractional operators in data science, physics,  and
machine learning.

\begin{figure}[h!]
\label{f:Gauss}
	\includegraphics[width=\textwidth]{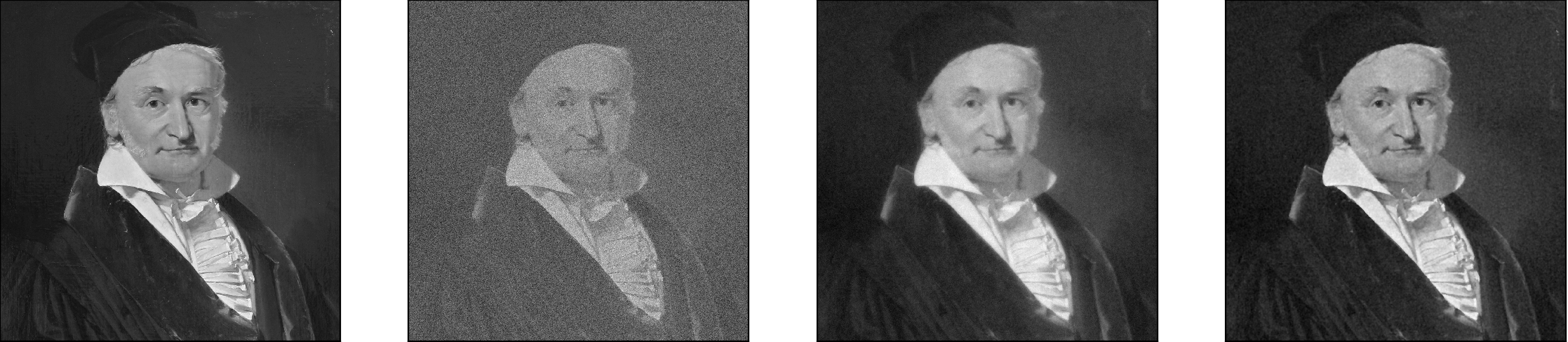}
	\caption{From left to right, original image, noisy image, denoised image using 
	fractional regularization with randomly selected regularization parameter and 
	denoised image using TV regularization with optimized regularization parameter.
	The figures have been reproduced from \cite{HAntil_SBartels_2017a}.}
\end{figure}

\medskip
\noindent 
{\bf Fractional Laplacian in imaging:}
%
Image denoising problems have received a lot of attention over last several decades. 
A major breakthrough occured when the article \cite{LIRudin_SOsher_EFatemi_1992a} 
considered the Total Variation (TV) regularization to capture the jumps in an image 
reconstruction, during the noise removal process. Ever since, this approach has 
attracted a lot of attention from researchers interested in both theory and numerical 
methods. However, it is a still quite a challenging problem because of the non-smooth 
nature of TV and the fact that when we formally write the Euler-Lagrange equations, 
they are nonlinear and degenerate. Recently in \cite{HAntil_SBartels_2017a}, the 
authors replaced the total variation semi-norm by fractional Laplacian regularization
with fractional Laplacian $(-\Delta)^s$ of order $s$ with constant $s \in (0,1)$. The fractional 
variational model is given by 
	\begin{equation}
		\label{eq:imag0}
		\min_{u} \frac12 \int_\Om | (-\Delta)^{\frac{s}{2}} u |^2 \dv{x} 
			+ \frac{{\lambda}}{2} \| u - f \|^2_{L^2(\Om)} 
	\end{equation}
where $f$ is the noisy image, $\lambda$ is the regularization parameter, and $\Om$ is the 
image domain. Moreover, periodic boundary conditions are assumed. 
The key advantage is that the Euler-Lagrange equation in the case of \eqref{eq:imag0} is a 
linear elliptic fractional partial differential equation (PDE) of the type 
	\begin{equation}
		\label{eq:imag1}
		(-\Delta)^s u + \lambda u = \lambda f  \quad \mbox{in } \Om 
	\end{equation}
which can be easily solved in the case of periodic boundary conditions by using Fourier series \cite{HAntil_SBartels_2017a}, or for general zero exterior conditions by using the method of bilinear forms. 	
From left to right, Figure~\ref{f:Gauss} shows the original image, noisy image, denoised 
image using fractional model with fractional exponent of $s = 0.42$ and $\lambda = 10$.
The rightmost image has been obtained using TV regularization approach 
\cite{AChambolle_2004a,MBurger_KFrick_SOsher_OScherzer_2007a}, 
where we have even optimized over $\lambda$. It is clear that the reconstruction quality
in the fractional case is comparable to the TV case even without optimizing over 
$\lambda$, however, it is significantly cheaper than TV. We also refer to \cite{HAntil_ZDi_RKhatri_2020a} 
for the extension of this work to tomographic reconstruction problems where the second
term in \eqref{eq:imag0} has been replaced by $\| K u - f \|^2_{L^2(\Om)}$, with $K$
denoting a linear operator. The article \cite{HAntil_ZDi_RKhatri_2020a} also discusses 
a bilevel optimization and machine learning based strategy to identify parameters such
as $\lambda$ and $s$. See also \cite{MDElia_JCDLosReyes_AMTrujillo_2021a} for another 
recent work on bilevel strategies for nonlocal problems. 
 
\begin{figure}[h!]
	\centering
	\includegraphics[width=0.28\textwidth]{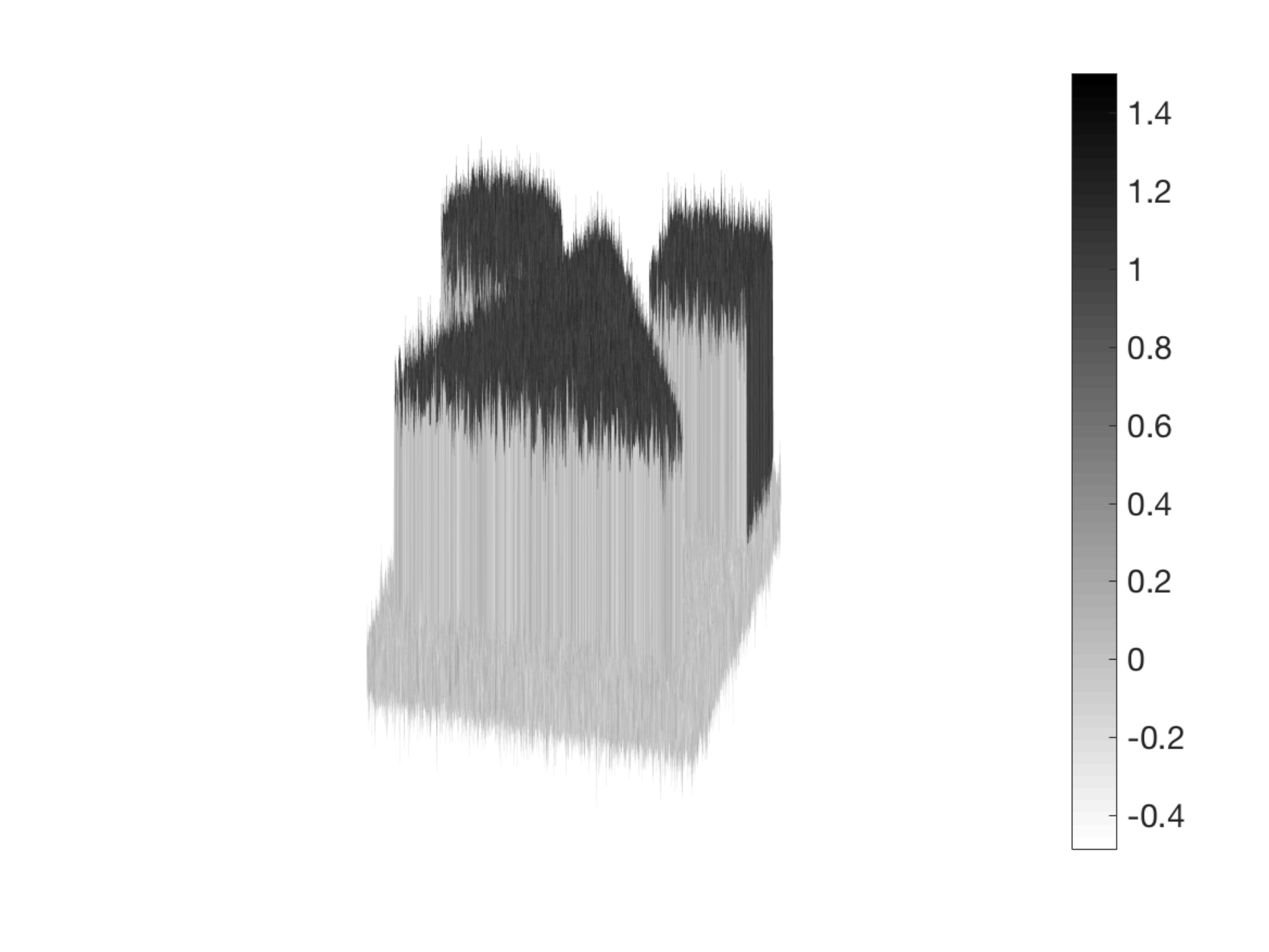}
	\includegraphics[width=0.28\textwidth]{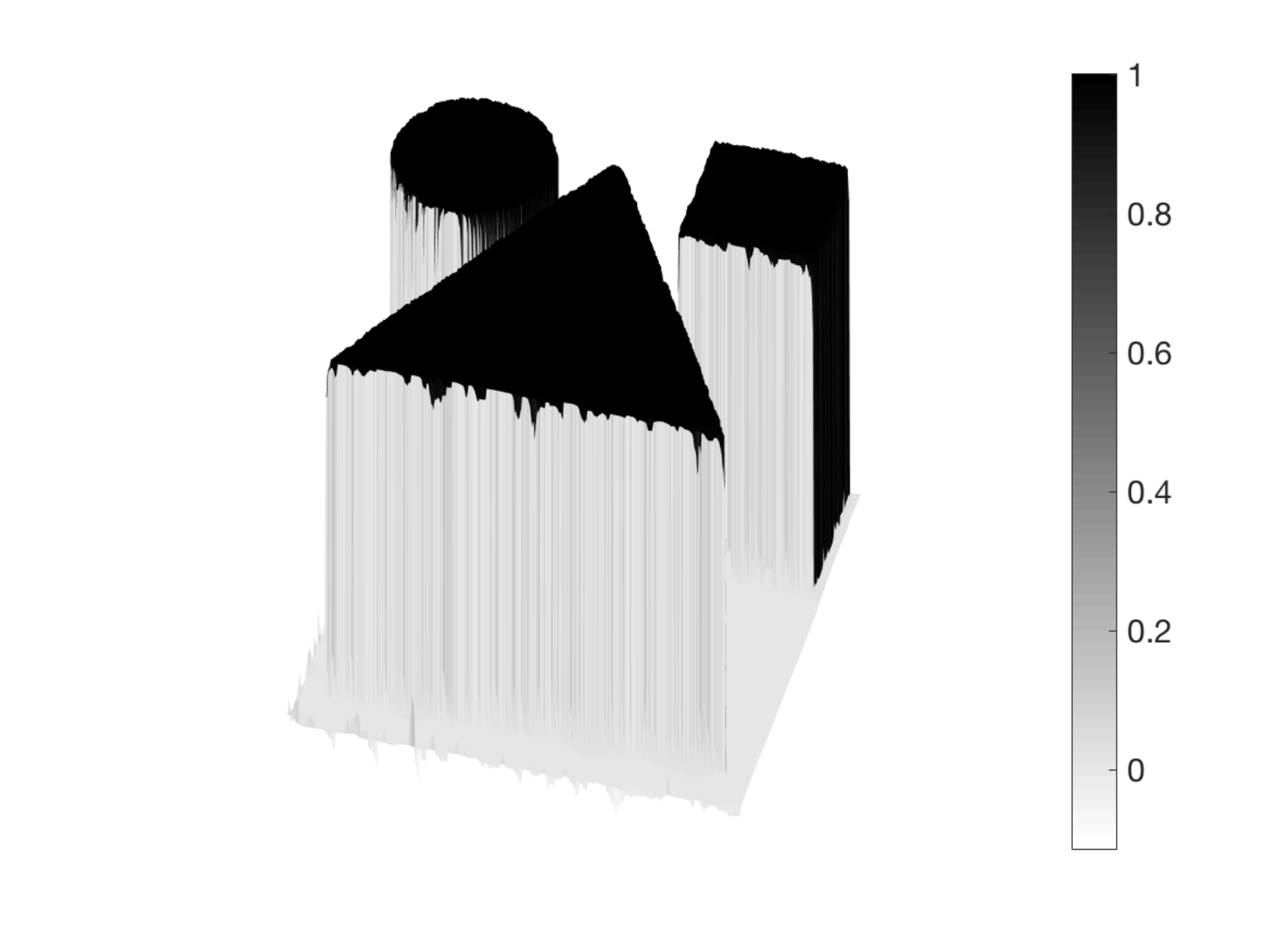}
	\includegraphics[width=0.28\textwidth]{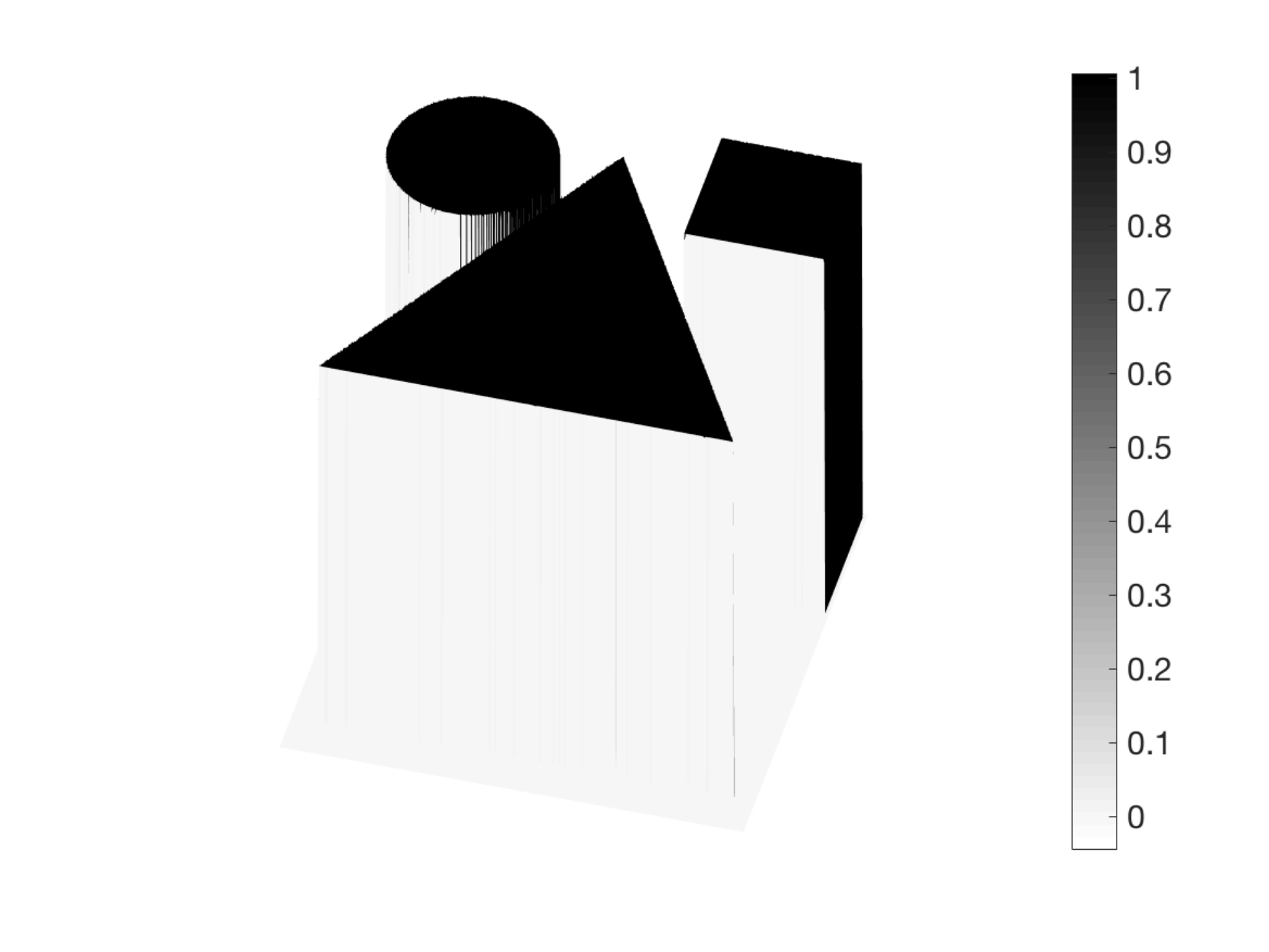}
	\caption{\emph{Left:} Noisy image. \emph{Middle:} reconstruction using TV regularization
	where the regularization parameter has been optimized. \emph{Right:} Reconstruction using $s(x)$ approach
	where the regularization parameter $\lambda$ has been chosen randomly. 
	Notice that the reconstruction using $s(x)$ approach is almost perfect. 
	The figures have been reproduced from \cite{HAntil_CNRautenberg_2019a}.}
	\label{f:sx}
\end{figure}

The article \cite{HAntil_CNRautenberg_2019a} introduces a new model with variable 
order $s(x) \in [0,1]$ with $x \in \O$, i.e., $s$ is not a constant any longer but instead, it is a function. 
Moreover, $s$ is allowed to take the extreme values of 0 and 1. This can be
quite beneficial for applications (e.g. imaging) where it is critical
to capture the jump across the interfaces (e.g., phase field models). Figure~\ref{f:sx} 
shows a comparison between TV regularization and the variable order approach.
A strategy to identify the function $s$ is also discussed in \cite{HAntil_CNRautenberg_2019a}. We 
clearly notice that TV is rounding up the corners in the middle panel, while the variable
$s$ based approach in the right panel leads to an almost perfect reconstruction. Notice that the regularization 
parameter has been optimized in TV case and has been chosen randomly in the variable $s$ 
case.

\begin{figure}[h!]
	 \centering
	\includegraphics[width=0.35\textwidth]{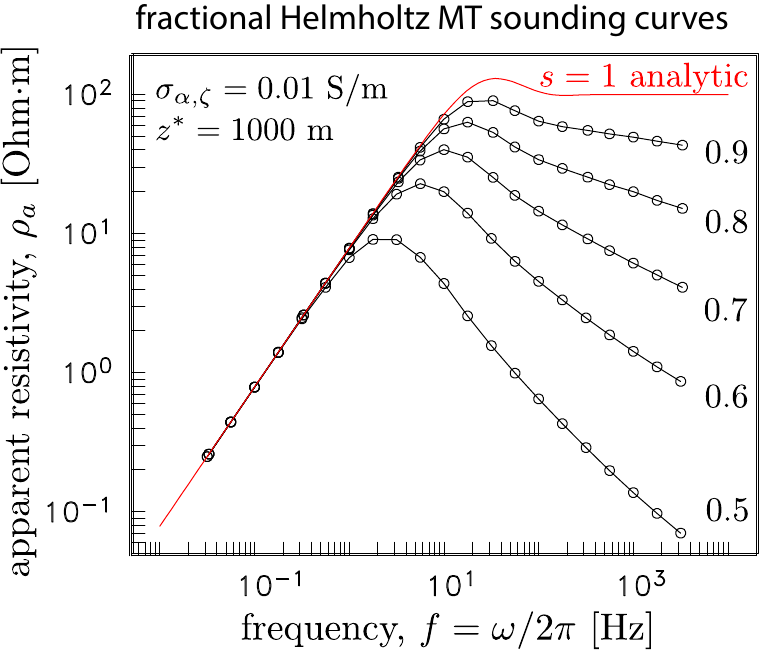}  \quad       
        \includegraphics[width=0.5\textwidth,height=.33\textwidth]{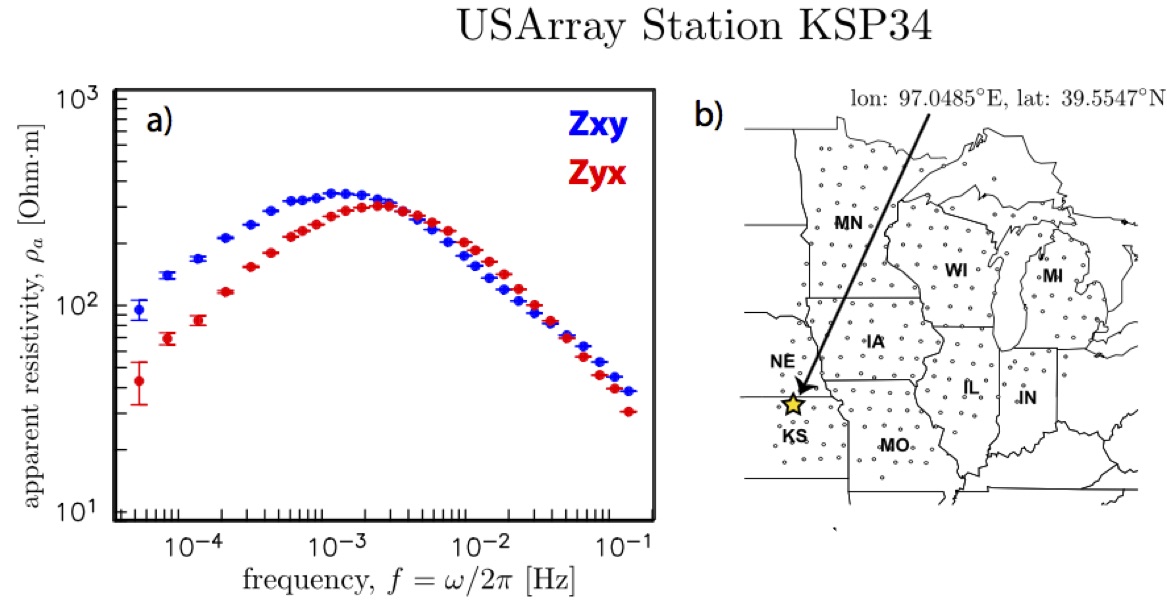} 
        \caption{\emph{Left:} Results from numerical simulations. \emph{Right:} The apparent resistivity data from 
        USArray MT station for KSP34 located NW of Kansas City, KS, USA from the US Array.}
        \label{f:fracHelm}
\end{figure}

\medskip
\noindent 
{\bf Fractional Laplacian in geophysics:} 
%
The article \cite{CJWeiss_BGVBWaanders_HAntil_2020a} has recently derived the fractional Helmholtz equation
using the first principle arguments in conjunction with a constitutive relationship. A finite element based numerical 
method motivated by \cite{ABonito_JEPasciak_2015a} has been introduced to solve the problem. These concepts 
are applied to the scalar Helmholtz equation and its use in electromagnetic interrogation of Earth's interior through 
the magnetotelluric method. For the magnetotelluric problem, several interesting features are  observable in the field data: long-range correlation of the predicted electromagnetic fields; a power-law relationship between the squared impedance amplitude and squared wavenumber whose slope is a function of the fractional 
exponent within the governing Helmholtz equation; and, a non-constant apparent resistivity spectrum whose 
variability arises solely from the fractional exponent. The latter can be seen in Figure~\ref{f:fracHelm} (left). The 
figure on the right is the apparent resistivity data from USArray MT station for KSP34 located NW of Kansas City, 
KS, USA from the US Array. Notice that the fractional model provides an excellent qualitative match.

Other applications of fractional Laplacian include, \emph{harmonic maps} \cite{HAntil_SBartels_ASchikorra_2021a}
and \emph{quantum spin chains} \cite{FDMHaldane_1986a} etc.

\medskip
\noindent 
{\bf Fractional time derivative in machine learning:} 
%
Dynamical system based deep neural networks (DNNs) are gaining popularity, as they allow us to model 
connectivity among the network layers. Recently in \cite{HAntil_RKhatri_RLohner_DVerma_2020a}, we 
have considered  
a new framework for classification problems and it has been further extended to parameterized 
partial differential equations in \cite{HAntil_HCElman_AOnwunta_DVerma_2021a}. 
These works are motivated by \cite{LRuthotto_EHaber_2019a}. These articles attempt to develop 
mathematical models for the analysis and understanding of DNNs. The idea is to consider DNNs 
as dynamical systems. In particular, in 
\cite{MBenning_ECelledoni_MEhrhardt_BOwren_CBSchnlieb,SGuenther_LRuthotto_JSchroder_ECyr_NGauger_2020a}, 
the process of training a 
DNN is thought of as an optimization problem constrained by a discrete ordinary differential 
equation. Designing DNN algorithms at the continuous level has the appealing 
advantage of architecture independence; in other words, the number of optimization 
iterations remains stable as the number of layers are increased. The articles 
\cite{HAntil_RKhatri_RLohner_DVerma_2020a,HAntil_HCElman_AOnwunta_DVerma_2021a}
specifically consider fractional ODE as constraints. This fractional DNN allows the 
network to access historic information of input and gradients across antecedent layers
since each layer is connected to all previous layers unlike the standard DNN, where 
each layer is connected only to a single previous layer.

\medskip
\noindent 
{\bf Outline:} In view of the above motivations, this article focuses on several optimization problems
constrained by fractional partial differential equations. We will use the term ``optimal control" to 
describe such problems. In Section~\ref{s:two} we begin by defining 
two fractional operators, i.e., fractional Laplacian and fractional time derivative. Section~\ref{s:fracdiff}
focuses on various notions of solutions to a linear elliptic fractional diffusion equation. We discuss 
various notions of solutions, including weak and very-weak solutions under various regularity conditions
on data. We also provide a discussion on various approaches to approximate the fractional Laplacian
and fractional time derivative. Section~\ref{s:ext} focuses on a new notion of optimal control, i.e., 
exterior optimal control. Moreover, Section~\ref{s:stateconst} focuses on elliptic fractional state 
constraint problems where we also discuss a Moreau-Yosida based algorithm to solve the problem. 
Section~\ref{s:time} is devoted to fractional deep neural networks which is also an optimal control
problem governed by time fractional ordinary differential equation.  
We conclude by stating various open problems in Section~\ref{s:open}. 
We emphasize that this is a review article and the results have been previously
published in the authors' other research papers.

\section{Two fractional operators and their properties}
\label{s:two}

Unless otherwise stated, throughout the paper, we assume that $\Om \subset \R^N$ is a bounded
open Lipschitz domain with boundary $\partial\Om$. Notice that many of our results hold without
the Lipschitz assumption, but we do not get into those details here in this review. 



\subsection{Fractional Laplacian $(-\Delta)^s$ and nonlocal normal derivative $\mathcal{N}_s$}

%

To define the integral fractional Laplace operator we consider the weighted Lebesgue space
\[
\mathbb{L}^1_s(\mathbb{R}^N) = 
\Big\{ f : \mathbb{R}^N \rightarrow \mathbb{R} \mbox{ measurable and } \int_{\mathbb{R}^N} \frac{|f(x)|}{(1+|x|)^{N+2s}} \dv{x} < \infty \Big\} ,
\]
and first define for $f \in \mathbb{L}^1_s(\mathbb{R}^N)$, $\varepsilon > 0$, and $x\in \mathbb{R}^N$ the quantity
\[
(-\Delta)^s_\varepsilon f (x) = C_{N,s} \int_{\{ y \in \mathbb{R}^N , |y-x| > \varepsilon \}} \frac{f(x)-f(y)}{|x-y|^{N+2s}} \dv{y} ,
\]
where the constant $C_{N,s} = \big(s2^{2s}\Gamma\left(\frac{2s+N}{2}\right)\big)/\big(\pi^{\frac{N}{2}}\Gamma(1-s)\big)$	
is obtained by Euler Gamma function. 
We then define the  {\em integral version} of the fractional Laplace operator for $s \in (0,1)$
via a limit passage for $\veps \to 0$, i.e., 
\begin{equation}\label{eq:integ}
(-\Delta)^s f(x) = C_{N,s} \mbox{P.V. } \int_{\mathbb{R}^N}  \frac{f(x)-f(y)}{|x-y|^{N+2s}} \dv{y} 
	       = \, \lim_{\varepsilon \to 0} (-\Delta)^s_{\varepsilon} f(x) ,
\end{equation}
where P.V. indicates the Cauchy principal value. Note that this definition for the full space $\mathbb{R}^N$ 
coincides with the spectral definition of the fractional Laplacian obtained using Fourier transform 
\cite[Proposition~3.4]{DDiNezza_GPalatucci_EValdinoci_2012a}, see also \cite{LCaffarelli_LSilvestre_2007a}. 
Such an equivalence also holds in the case of periodic boundary conditions 
\cite[Eq.~(2.53)]{NAbatangelo_EValdinoci_2019a}. 

\medskip
\noindent
{\bf Fractional order Sobolev spaces:} 
Next we introduce the fractional order  
Sobolev space $H^s(\mathbb{R}^N)$ for $s \in (0,1)$ by setting 
\begin{align*}
	[f]_{H^{s}(\R^N)}   &= \| (-\Delta)^{\frac{s}{2}} f\|_{L^2(\R^N)} = \left( \int_{\R^N} \int_{\R^N} \frac{|f(x)-f(y)|^2}{|x-y|^{N+2s}} \dv{y} \dv{x} \right)^{\frac12},\\
	\|f\|_{H^{s}(\R^N)} &= \|f\|_{L^2(\R^N)} + [f]_{H^{s}(\R^N)} .
\end{align*}
Then, the Sobolev space $H^{s}(\R^N)$ defined as 
\[
	H^{s}(\R^N) = \left\{  f \in L^2(\R^N) \, : \, \|f\|_{H^{s}(\R^N)} < +\infty \right\}
\]
is a Hilbert space. We denote the dual of $H^s(\R^N)$ by $H^{-s}(\R^N)$ and
$\langle \cdot, \cdot \rangle$ the duality pairing between $H^s(\R^N)$ and 
$H^{-s}(\R^N)$. 

For bounded open sets $\O\subset \R^N$ and parameters $s \in (0,1)$ we define
Sobolev spaces $\widetilde{H}^{s}(\Om)$ by considering trivial extensions to $\R^N$, i.e.,
we set  
\[ 
\widetilde{H}^{s}(\Om) = \{f \in L^2(\R^N) \, : \, (-\Delta)^{\frac{s}{2}} f \in L^2(\R^N) , 
\quad f \equiv 0 \mbox{ in } \R^N \setminus \Om  \} .
\]
We recall the following density result for $s \in (0,1]$ and bounded domains $\Om \subset \R^d$ with a
Lipschitz continuous boundary \cite{AFiscella_FServadei_EValdinoci_2015a}
\[
\widetilde{H}^{s}(\Om) = \overline{\mathcal{D}(\Om)}^{\|\cdot\|_{\widetilde{H}^{s}(\Om)}} ,
\]
and by a Poincar\'e type inequality, which is a consequence of H\"older's inequality and Sobolev 
imbedding theorem, \cite[Theorem 3.1.4.]{DRAdams_LIHedberg_1996a}, a norm on 
$\widetilde{H}^s(\Omega)$ is given by
\[
\|f\|_{\widetilde{H}^s(\Omega)} = \| (-\Delta)^{\frac{s}{2}} f\|_{L^2(\R^N)} .
\]

\medskip
\noindent
{\bf Other fractional operators:} 
Next we define the operator $(-\Delta)_D^s$ in $L^2(\Omega)$ as follows:
\begin{equation}
\label{eq:DLap}
	D((-\Delta)_D^s):=\Big\{u|_{\Om},\; u\in \widetilde{H}^{s}(\Om):\; (-\Delta)^su\in L^2(\Omega)\Big\},\;\;(-\Delta)_D^s(u|_{\Om}):=(-\Delta)^su\;\mbox{ a.e. in }\;\Omega.
\end{equation}
Notice that $(-\Delta)_D^s$ is the realization of the fractional Laplace operator $(-\Delta)^s$ in $L^2(\Om)$ 
with the zero Dirichlet exterior condition $u=0$ in $\RR^N\setminus\Omega$. 


For $u \in H^s(\R^N)$, we define the \emph{nonlocal normal derivative} $\mathcal{N}_s$ as: 
	\[
		\mathcal{N}_s u(x) := C_{N,s} \int_\O \frac{u(x) - u(y)}{|x-y|^{N+2s}} \dv{y}, \quad x \in \R^N \setminus \overline\O ,
	\]
which maps continuously $H^s(\R^N)$ into $H^{s}_{\rm loc} (\R^N\setminus \O)$. As a result, if $u \in H^s(\R^N)$, then 
$\mathcal{N}_s u \in L^2(\R^N\setminus \O)$. See \cite[Lemma~2.1]{HAntil_RKhatri_MWarma_2019a} for details. 

Moreover, the following \emph{integration-by-parts formula} can be found in \cite[Proposition~2.2]{HAntil_RKhatri_MWarma_2019a}. 
If $u \in H^s(\R^N)$ is such that $(-\Delta)^s u \in L^2(\O)$,  then for every $v \in H^s(\R^N)$ we have that 
	\begin{align}\label{Int-Part}
		\frac{C_{N,s}}{2} 
		\iint_{\RR^{2N}\setminus(\RR^N\setminus\Om)^2} 
		\frac{(u(x)-u(y))(v(x)-v(y))}{|x-y|^{N+2s}} \;dxdy 
		= \int_\Om v(-\Delta)^s u\;dx + \int_{\RR^N\setminus\Om} v\mathcal{N}_s u\;dx ,
	\end{align}
	where $\RR^{2N}\setminus(\RR^N\setminus\Om)^2 := (\Om\times\Om)\cup(\Om\times(\RR^N\setminus\Om))\cup((\RR^N\setminus\Om)\times\Om)$.

\subsection{Fractional time derivative: $_C\partial_t^\gamma$}

The focus of this section is to define the (strong) fractional order Caputo derivative and to introduce
the integration-by-parts formula. The latter is critical to drive the optimality conditions. 

\begin{definition}[Strong Caputo fractional derivative] 
\label{def:Caputo}
Let $u \in W^{1,1}([0,T];X)$, with $X$ denoting a Banach space. The (strong) Caputo fractional 
derivative of order $\gamma \in (0,1)$ is given by
	\begin{equation} \label{Lcaputo}
		_C\partial^{\gamma}_{t} u(t) 
		=\frac{1}{\Gamma(1-\gamma)} \int_{0}^{t}\frac{u'(r)}{(t-r)^{\gamma}} \dv{r} .
	\end{equation}
\end{definition}
The following integration-by-parts formula holds, \cite{OPAgrawal_2007a}, under the assumptions
of Definition~\ref{def:Caputo}
	\begin{equation}\label{IPF}
		\int_0^T v(t) _C\partial_t^\gamma u(t) \dv{t} = 
		\int_0^T \partial_{t,T}^\gamma v(t) u(t) \dv{t} + \left[ (I^{1-\gamma}_{t,T}v)(t) u(t) \right]_{t=0}^{t=T} ,
	\end{equation}
provided the expressions on the left and right hand sides make sense.  Here, 
$\partial_{t,T}^\gamma$ is the right Riemann-Liouville fractional derivative 
	\[
		\partial_{t,T}^\gamma u(t) = -\frac{d}{dt} (I^{1-\gamma}_{t,T}u)(t), \quad 
		\mbox{where} \quad
		I^{\gamma}_{t,T} u(t) := \frac{1}{\Gamma(\gamma)} \int_t^T (r - t)^{\gamma-1} u(r) \dv{r} 
	\]
denotes the right Riemann-Liouville fractional integral of order $\gamma$. We refer to 
\cite[pg.~76]{SGSamko_AAKilbas_OIMarchev_1993a} for further details.  We notice that if the space $X$ has the Radon-Nikodym property, then $W^{1,1}([0,T];X)$ is the optimal space for which \eqref{IPF} makes sense (see e.g.  the monograph \cite{GW-book} for more details).

\section{Fractional diffusion equation: analysis and numerical approximation}
\label{s:fracdiff}

\subsection{Non-homogeneous diffusion equations: analytic results}

The focus of this section is on stating well-posedness results for the linear elliptic fractional equation
\begin{equation}\label{eq:ell_d_gen}
	\begin{cases}
		(-\Delta)^s u = f \quad &\mbox{in } \Om, \\
		u = g \quad &\mbox{in } \RR^N\setminus \Om ,
	\end{cases}                
\end{equation} 
under various general assumptions on the data $f$ and $g$. In particular, we first establish the notion
of weak solution in a Sobolev space when $g$ is an appropriate Sobolev function itself. When 
$g$ is only square integrable, 
we establish the notion of very weak solution in $L^2$. Next for $g \equiv 0$, we provide conditions
on $f$ and $\O$ which leads to boundedness and continuity of $u$. This is followed by the notion of
very weak solution when $f$ is a Radon measure. We will denote the space of Radon measures
by $\mathcal{M}(\O)$. For parabolic versions of these results, we refer to \cite{HAntil_DVerma_MWarma_2020a,HAntil_DVerma_MWarma_2020b}. 
We further refer to \cite{HAntil_MWarma_2020a} for semilinear problems, 
\cite{HAntil_MWarma_2019a} for quasi-linear problems, \cite{HAntil_CNRautenberg_2018a,HAntil_CNRautenberg_ASchikorra_2020a} for variational and 
quasi-variational inequalities. 

We begin by stating the various notions of solutions to \eqref{eq:ell_d_gen}. This will be followed
by several well-posedness results. 
\begin{definition}[\bf Solutions to elliptic Dirichlet problem] 
	\label{def:weak_ell_d}
	
We define various notions of solutions to \eqref{eq:ell_d_gen}: 		
\begin{enumerate}[(i)]	
	\item {\bf Weak solution to non-homogeneous problem:} 
	Let $f \in \widetilde{H}^{-s}(\Om)$, $g \in H^s(\RR^N\setminus\Om)$ and 
	$\mathcal{Z} \in H^s(\RR^N)$ be such that $\mathcal{Z}|_{\RR^N\setminus\Om} = g$. 
	A function $u \in H^s(\RR^N)$ is said to be a weak solution to \eqref{eq:ell_d_gen} if 
	$u-\mathcal{Z} \in \widetilde{H}^{s}(\Om)$ and 
	\begin{equation}\label{eq:weak}
	\frac{C_{N,s}}{2} 
	\int_{\RR^N}\int_{\RR^{N}} 
	\frac{(u(x)-u(y))(v(x)-v(y))}{|x-y|^{N+2s}} \dv{x} \dv{y} 
	= \langle f,v \rangle  , \quad \forall v \in \widetilde{H}^{s}(\Om) . 
	\end{equation}
	
	\item {\bf Very-weak solution to non-homogeneous Dirichlet problem:} 
	Let $g \in L^2(\RR^N\setminus\Om)$ and $f \in \widetilde{H}^{-s}(\Om)$. 
	A function $u \in L^2(\RR^N)$ is said to be a very-weak solution to \eqref{eq:ell_d_gen} if the identity
	\begin{equation}\label{eq:vw_ell_d}
		\int_{\Om} 
		u (-\Delta)^s v\;dx
		= \langle f,v \rangle - \int_{\RR^N\setminus\Om} g \mathcal{N}_s v\;dx,
	\end{equation}
	holds for every $v \in V := \{v \in \widetilde{H}^{s}(\Om) \;:\; (-\Delta)^s v
	\in L^2(\Om) \}$. 
	
	\item {\bf Very-weak solution to homogeneous Dirichlet problem with measure data:} 
	Let $p$ satisfy  
	\begin{equation}\label{cond-p}
			p>\frac{N}{2s}\;\;\mbox{ if }\; N>2s,\quad
			p>1 \;\;\mbox{ if }\; N=2s,\quad
			p=1\;\;\mbox{ if }\; N<2s,
		\end{equation}  
	 and $\frac{1}{p}+\frac{1}{p'} = 1$. 
	Let $f \in \mathcal{M}(\Omega)$.
	A function $u \in L^{p'}(\Om)$ is said to be a very-weak solution 
	to \eqref{eq:ell_d_gen}, if  for every $v \in V := \{v \in C_0(\Om) \cap \widetilde{H}^{s}(\Om):\; (-\Delta)^s v \in L^p(\Om) \}$ we have
	\begin{equation}\label{EWW}
		\int_{\Om} 
		u(-\Delta)^s v\;dx
		= \int_\Om v \; df .
	\end{equation}
Here $C_0(\Om)$ is the space of all continuous functions in $\overline\Om$ that vanish
on $\pOm$. 
\end{enumerate}	
\end{definition}

%

\begin{theorem}
	\label{thm:not_soln_ellip}
	The following results hold for non-homogeneous Dirichlet boundary value problems:
	\begin{enumerate}[(i)]
		\item {\bf Weak solution:} Given $f \in \widetilde{H}^{-s}(\O)$ and $g \in H^{s}(\RR^N \setminus\Om)$
			 there exists a unique weak solution to \eqref{eq:ell_d_gen} according to Definition~\ref{eq:ell_d_gen} (i).
			Also there is a constant $C > 0$ such that 
			\begin{align}\label{Es-DS_ell}
				\|u\|_{H^s(\RR^N)} 
				\le C \left(\|f\|_{\widetilde{H}^{-s}(\Om)} 
				+ \|g\|_{H^s(\RR^N \setminus\Om)} \right).
			\end{align}
		
		\item {\bf Very-weak solution:} Given $f \in \widetilde{H}^{-s}(\O)$, $g \in L^2(\R^N\setminus \O)$
			then there exists a unique very-weak solution $u$ to 
			\eqref{eq:ell_d_gen} according to Definition~\ref{eq:ell_d_gen} (ii) that fulfills
			\begin{align}\label{VWS_EST_ell}
				\|u\|_{L^2(\Om)} \le C\left(\|f\|_{\widetilde{H}^{-s}(\Om)} 
				+ \|g\|_{L^2(\RR^N \setminus\Om)} \right),
			\end{align}
			and the following assertions hold: 
			(a) Every weak solution of \eqref{eq:ell_d_gen} is also a very-weak solution.
			(b) Every very-weak solution of \eqref{eq:ell_d_gen} that belongs to $H^s(\RR^N)$ is also a weak solution.			
	\end{enumerate}
	In addition, the following holds true when $g \equiv 0$:
	\begin{enumerate}[(i)]
		\item[(iii)] {\bf $u$ is bounded in $L^\infty(\Om)$:} Let $f \in L^p(\Om)$ with $p$ as in \eqref{cond-p}, 
		then $u \in L^\infty(\Om) \cap \widetilde{H}^s(\Om)$. Moreover, if $\Om$ fulfills the exterior cone condition, then the weak
		solution $u \in C_0(\Om)$ and there is a $C = C(N,s,p,\Om) > 0$ such that $\| u \|_{C_0(\Om)} \le C \|f\|_{L^p(\Om)}$. 
		
		\item[(iv)] {\bf Very-weak solution with measure $f$:} Let $f \in \mathcal{M}(\Om)$, $p$ as in \eqref{cond-p}, 
		and $\O$ fulfills the exterior cone condition, then there exists a unique $u \in L^{p'}(\Om)$ solving \eqref{eq:ell_d_gen} 
		according to Definition~\ref{eq:ell_d_gen} (iii), and there is a constant $ C = C(N,s,p,\Om) > 0$ such that $\|u\|_{L^{p'}(\O)} 
		\le C \|f\|_{\mathcal{M}(\Om)}$. 
	\end{enumerate}	
\end{theorem}

\subsection{Non-homogeneous diffusion equations: numerical approximation}
\label{s:nhbc}

Approximation of fractional Laplacian is a challenging topic and has received a tremendous
amount of attention recently. However, most of the focus has been on the homogeneous problem. 
There exist several efficient works in 1D but for $N > 1$ 
the problem is very challenging and the number of works are limited. The main difficulty 
is the fact that the bilinear form \eqref{eq:weak} 
contains a singular integral and the traditional finite element approaches are not 
amenable to this. The first work that rigorously provides approximation to the weak
solutions is by  Acosta and Borthagaray \cite{GAcosta_JPBorthagaray_2017a}, see also 
\cite{GAcosta_FMBersetche_JPBorthagaray_2017a}. 
However, implementations in these works have been limited to $N=2$ dimensions. 
We also refer to another finite element approach by Bonito, Lei, and Pasciak 
\cite{ABonito_WLei_JEPasciak_2019a}, which also works for $N = 3$. 

In contrast, in \cite{HAntil_PDondl_LStriet_2021a} the authors have introduced an efficient  
spectral method which works in arbitrary Lipschitz domains. Under this method, the evaluation 
of the fractional Laplacian and its application onto a vector has complexity of 
$\mathcal{O}(M\log(M))$ where $M$ is the number of unknowns. For instance, for exponent 
$s = 1/4$ the 3D implementation can solve the Dirichlet problem with $5 \cdot 10^6$ unknowns 
under 2 hours on a standard office workstation. 

The spectral method presented in \cite{HAntil_PDondl_LStriet_2021a} also extends to
the non-homogeneous Dirichlet problem. But below, instead, we briefly discuss a finite
element method based on \cite{HAntil_RKhatri_MWarma_2019a,HAntil_DVerma_MWarma_2020a}. 
We refer to \cite{GAcosta_JPBorthagaray_NHeuer_2019a} for an alternative approach. 
The main idea in \cite{HAntil_RKhatri_MWarma_2019a,HAntil_DVerma_MWarma_2020a} is to 
approximate the solutions to a Dirichlet problem by solutions to a Robin problem. Some of the
key challenges for the Dirichlet problem are:
%
	\begin{enumerate}[$\bullet$]
		\item Since we are interested in optimal control problems with $L^2$-exterior controls,
			 the correct notion of solution is as given in \eqref{eq:vw_ell_d}. Numerically,
			this will require approximating $\mathcal{N}_s$.
			
		\item The control equation, in addition, will require approximation of $\mathcal{N}_s$ 
			applied to the adjoint variable. 
	\end{enumerate} 
In other words, we have to approximate $\mathcal{N}_s$ in addition to approximating $(-\Delta)^s$. 
Approximating the Dirichlet problem by the Robin problem helps overcome both these issues. To define the Robin
problem, we consider the Sobolev space introduced in 
\cite{SDipierro_XRosOton_EValdinoci_2017a}. Let $g\in L^1(\RR^N\setminus\Omega)$ be fixed and set 
\begin{align*}
	W_{\Omega,g}^{s,2}:=\Big\{u:\RR^N\to\RR\;\mbox{ measurable and }\,\|u\|_{W_{\Omega,g}^{s,2}}<\infty\Big\}, 
\end{align*}
where for the measure $\mu$ on $\R^N \setminus \O$ given by $\dv{\mu} = |g| \dv{x}$, we have 
\begin{align}\label{norm-e}
	\|u\|_{W_{\Omega,g}^{s,2}}:=\left(\|u\|_{L^2(\Omega)}^2+\|u\|_{L^2(\RR^N\setminus\Omega,\mu)}^2+\int\int_{\RR^{2N}\setminus(\RR^N\setminus\Omega)^2}\frac{|u(x)-u(y)|^2}{|x-y|^{N+2s}}dxdy\right)^{\frac 12}.
\end{align}
The article \cite[Proposition 3.1]{SDipierro_XRosOton_EValdinoci_2017a} shows that for 
$g\in L^1(\RR^N\setminus\Omega)$, $W_{\Omega,g}^{s,2}$ is a Hilbert space. Now for 
every $n \in \mathbb{N}$, we  define the generalized Robin problem  
\begin{equation}\label{eq:Sn2}
\begin{cases}
(-\Delta)^su=f\;\;&\mbox{ in }\;\Omega,\\
\mathcal N_su+ n \kappa u=n \kappa z\;\;\;&\mbox{ in } \;\RR^N\setminus\Omega ,
\end{cases}
\end{equation}
and the following result holds (see Proposition~\cite[3.9]{HAntil_RKhatri_MWarma_2019a})
\begin{proposition}
	Let $n \in \mathbb{N}$ and $\kappa\in L^1(\Omc)\cap L^\infty(\Omc)$.
	Then for every $z\in L^2(\RR^N\setminus\Omega,\mu)$ and $f\in ( W_{\Om,\kappa}^{s,2})^\star$, 
	there exists a  weak solution $u\in W_{\Om,\kappa}^{s,2}$ of \eqref{eq:Sn2} in the following sense:
	\begin{align}\label{we-so}
		\frac{C_{N,s}}{2}\int\int_{\RR^{2N}\setminus(\Omc)^2}&\frac{(u(x)-u(y))(v(x)-v(y))}{|x-y|^{N+2s}}\;dxdy+n \int_{\Omc}\kappa uv\;dx\notag\\
		&=\langle f,v\rangle_{ ( W_{\Om,\kappa}^{s,2})^\star, W_{\Om,\kappa}^{s,2}}+ n \int_{\Omc}\kappa zv\;dx.
	\end{align}
\end{proposition}
Next we make the following assumption.
\begin{assumption}\label{asum}
We assume that $\kappa\in L^1(\Omc)\cap L^\infty(\Omc)$ and satisfies $\kappa>0$ almost everywhere in $K:=\mbox{supp}[\kappa]\subset\Omc$, where the Lebesgue measure $|K|>0$.
\end{assumption}
Under this assumption, the solution to \eqref{we-so} belongs to $W_{\Om,\kappa}^{s,2} \cap L^2(\R^N\setminus\O)$ \cite[Lemma~6.2]{HAntil_RKhatri_MWarma_2019a}.
Moreover, the following approximation result holds
	\begin{theorem}[\bf Approximation of Dirichlet solution by Robin solution]
	\label{thm:dbc_approx}
		Under Assumption~\ref{asum}, the following assertions hold:
		\begin{enumerate}[(a)]
			\item  Let $z \in H^s(\Omc)$ and  $u_{n} \in W^{s,2}_{\Om,\kappa}\cap L^2(\Omc)$ be the weak solution 
				 of  \eqref{we-so}. Let $u\in H^s(\RR^N)$ be the weak solution to the state equation 
				 \eqref{eq:weak}. Then there is a constant $C>0$ (independent of $n$) such that
		                 \begin{align}\label{es-diff}
				  \|u-u_{n}\|_{L^2(\RR^N)}\le \frac{C}{n}\|u\|_{H^s(\RR^N)}.
			        \end{align}
			       In particular $u_{n}$ converges strongly to $u$ in $L^2(\RR^N)$ as $n\to\infty$.
			\item Let $z \in L^2(\Omc)$ and $u_{n} \in W^{s,2}_{\Om,\kappa}\cap L^2(\Omc)$ be the weak solution of \eqref{we-so}. 
				Then there exist a subsequence that we still denote by $\{u_n\}_{n\in\NN}$  and a function $\tilde u\in L^2(\RR^N)$ such that 
				$u_{n}\rightharpoonup \tilde u$ in $L^2(\RR^N)$ as $n\to\infty$, and $\tilde u$ satisfies
				\begin{align}\label{eq63}
					\int_{\Omega}\tilde u(-\Delta)^sv\;dx=-\int_{\Omc}\tilde u\mathcal N_sv\;dx,
				\end{align}
				for all $v\in V := \{v \in \widetilde{H}^{s}(\Om) \;:\; (-\Delta)^s v \in L^2(\Om) \}$.
		\end{enumerate}		
	\end{theorem}
Next, we introduce a discrete scheme to approximate \eqref{we-so}.  Let $\widetilde\O$ be an open 
bounded set that contains $\O$. We consider a conforming simplicial triangulation of $\O$ and 
$\widetilde\O \setminus \O$ such that the partition remains admissible. We assume that the support 
of $z$ and $\kappa$ is contained in $\widetilde\O \setminus \O$. We let the finite element space 
$\mathbb{V}_h$ on $\widetilde\O$ to be the set of continuous piecewise linear functions. Then it 
remains to approximate the weak form \eqref{we-so}. All other terms are standard and can be
done using any standard finite element code and quadrature rule, except the stiffness matrix. We 
assemble the latter using \cite{GAcosta_FMBersetche_JPBorthagaray_2017a}. All other matrices
are computed by using quadrature which is accurate for polynomials of degree less than or equal to 4.

\begin{figure}[htb]
  \centering
  \includegraphics[width=0.33\textwidth]{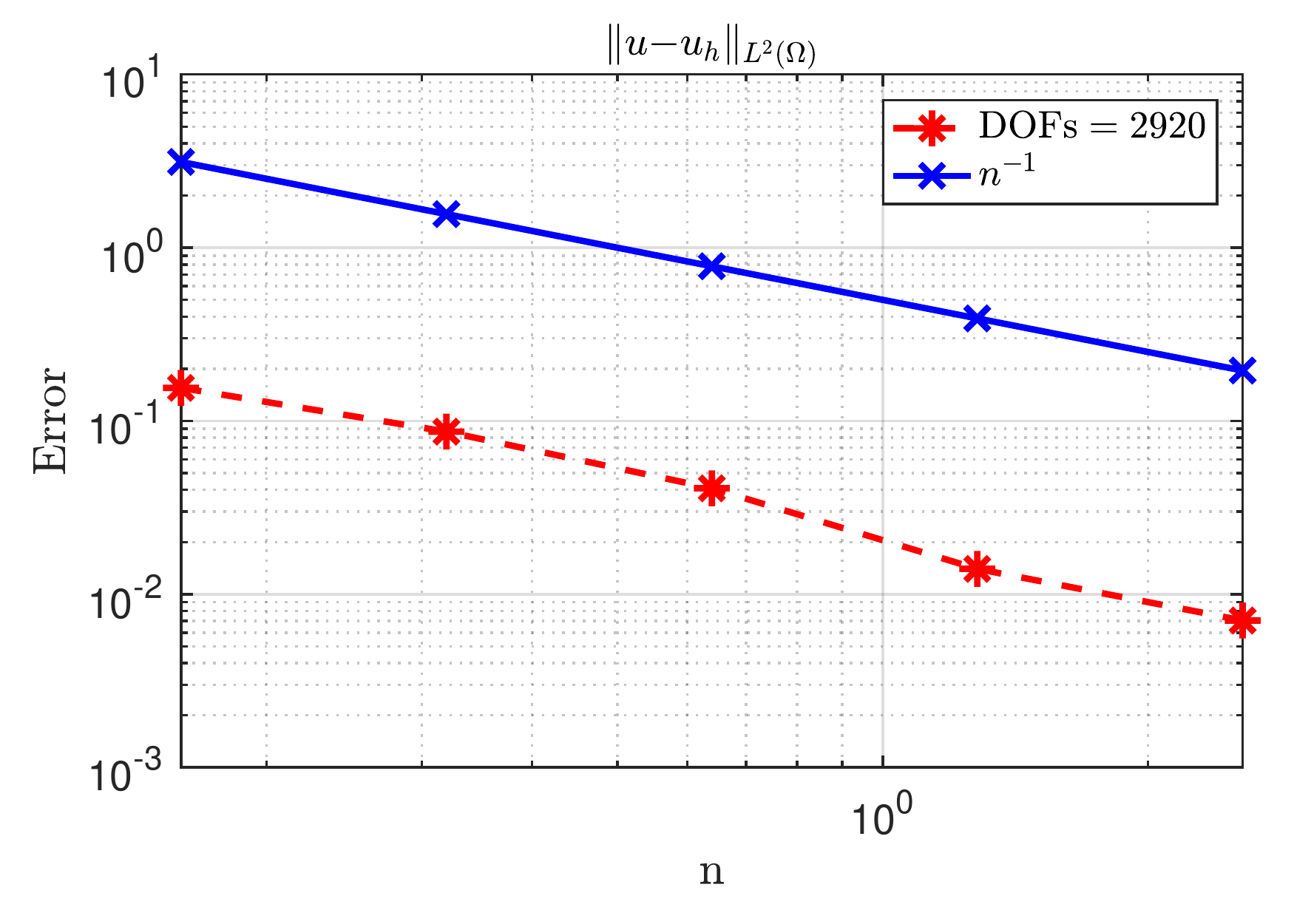}  
  \includegraphics[width=0.33\textwidth]{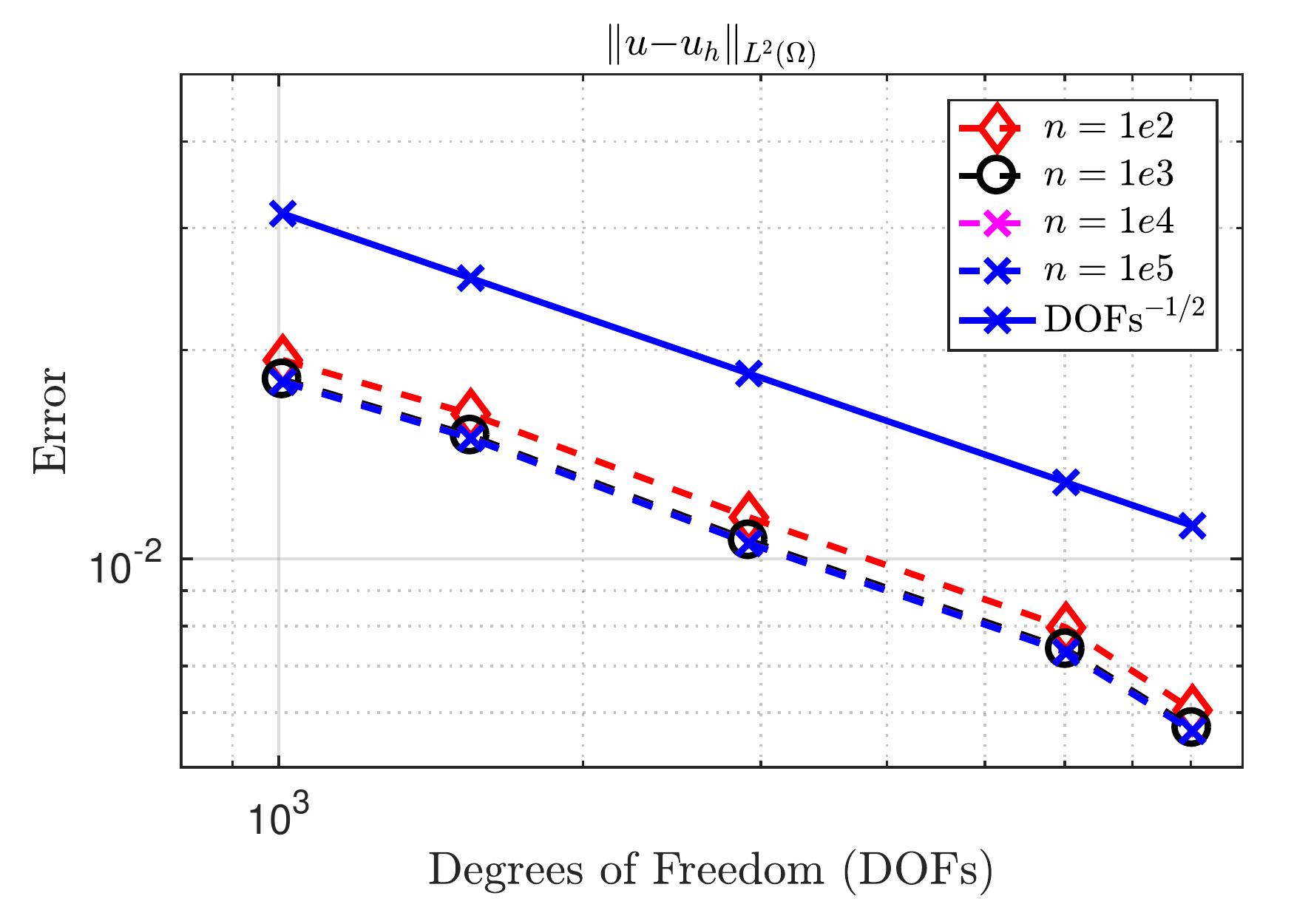} 
  \caption{\label{eq:fig_rate}
   Left panel: Let $s = 0.5$ and $\mbox{DoFs} = 2920$ be fixed. We let $\kappa = 1$ 
   and consider the $L^2$-error between the actual solution $u$ to the Dirichlet problem
   and its approximation $u_h$ which solves the Robin problem. We have plotted the 
   error with respect to $n$. The solid line denotes a reference line and the actual  
   error. 
   We observe a rate of $1/n$ which confirms our 
   theoretical result \eqref{es-diff}. Right panel: Let $s=0.5$ be fixed. For 
   each $n = 1e2, 1e3, 1e4, 1e5$ we have plotted the $L^2$-error with respect to 
   the degrees of freedom (DOFs) as we refine the mesh. Notice that the error is stable 
   with respect to $n$. Moreover, the observed rate of convergence is 
   $(\mbox{DoFs})^{-\frac12}$ and is independent of $n$.
   }
 \end{figure}

We consider an example taken from \cite{GAcosta_JPBorthagaray_NHeuer_2019a}. Let 
$\O = B_0(1/2) \subset \R^2$, i.e., a ball centered at 0 with radius 1/2 and $\widetilde\Om = B_0(3/2)$. 
Our goal is to find 
$u$ solving $(-\Delta)^s u = 2$ in $\O$ and $u(\cdot) = \frac{2^{-2s}}{\Gamma(1+s)^2} \left(1-|\cdot| \right)^s_+$ 
in $\R^N \setminus \O$. Figure~\ref{eq:fig_rate} shows our results and confirms our theoretical
findings in Theorem~\ref{thm:dbc_approx}.

\subsection{Fractional time derivative: numerical approximation}

We consider the $L^1$-scheme to approximate $_C\partial_t$ and we state the result 
with the help of a nonlinear ODE
\begin{equation} \label{eq:eq1}
	d_t^{\gamma}u(t)=f(u(t)), \quad u(0)=u_0.
\end{equation}
Consider the time discretization of the interval $[0,T]$ uniformly with step size $\tau$. We let 
$0 = t_0 < t_1 < \dots < t_{j+1} < \dots < t_N = T$, where $t_j = j \tau$. 

\begin{figure}[htb] 
		\centering
			\includegraphics[width=0.3\textwidth]{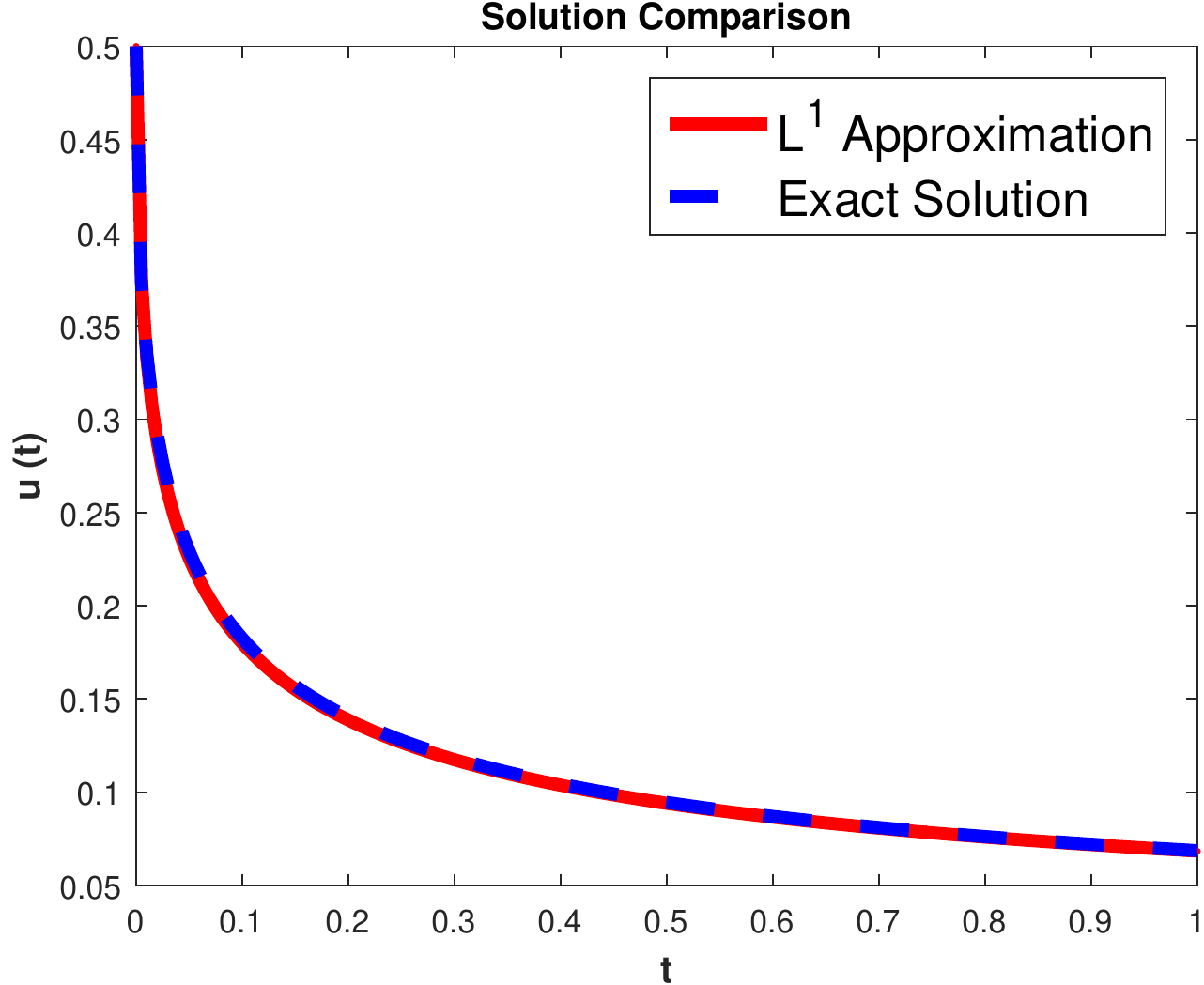}
			\caption{Comparison of exact solution of \eqref{ex:mittag} and its approximation using $L^1$-scheme.}
			\label{fig:GL_approx} 		
	\end{figure}

Then, using the $L^1$-scheme, the discretization of \eqref{eq:eq1} is given by: for $j = 0,...,N-1$,
\begin{alignat}{3} \label{eq:disc_cCaputo}
u(t_{j+1})= u(t_{j}) - \sum_{k=0}^{j-1}   a_{j-k}\: \left(u(t_{k+1}) - u(t_k) \right)  + \tau^{\gamma}\Gamma(2-\gamma)f(u(t_{j})),  
\end{alignat}
where the coefficients $a_{j-k}$ are given by, 
\begin{equation}\label{a_k}
a_{j-k}=(j+1-k)^{1-\gamma}-(j-k)^{1-\gamma}.
\end{equation}
Next, we present an example illustrating this approach. Consider the differential equation:
	\begin{align} \label{ex:mittag}
		d_t^{\frac12} u(t)=-4 u(t),\;\; u(0)=0.5.
	\end{align} 
Then, the exact solution to \eqref{ex:mittag} is given by, see \cite[Section 42]{SGSamko_AAKilbas_OIMarchev_1993a}, also 
\cite[Section 1.2]{IPodlubny1999a}, 
	$	u(t)=0.5\: E_{0.5}(-4t^{0.5}),$ 
where $E_{\alpha}$, with $0<\alpha \in \RR$, is the Mittag Leffler function, see \cite[Pg.~17]{IPodlubny1999a}, defined by
	$ 
		E_{\alpha }(z)=E_{\alpha , 1}(z)= \sum_0^{\infty} \frac{z^k}{\Gamma(\alpha k +1)} .
	$
Figure~\ref{fig:GL_approx} depicts the true solution and the numerical solution using the discretization \eqref{eq:disc_cCaputo} 
for the above example with uniform step size $\tau=0.005$ and final time, $T=1$. 
For numerical analysis of partial differential equation with strong Caputo derivatives, see  
\cite{RHNochetto_EOtarola_AJSalgado_2014b,HAntil_EOtarola_AJSalgado_2015a,BJin_RLazarov_ZZhou_2019a} and the references therein.


\section{Exterior optimal control of fractional parabolic PDEs with control constraints}
\label{s:ext}

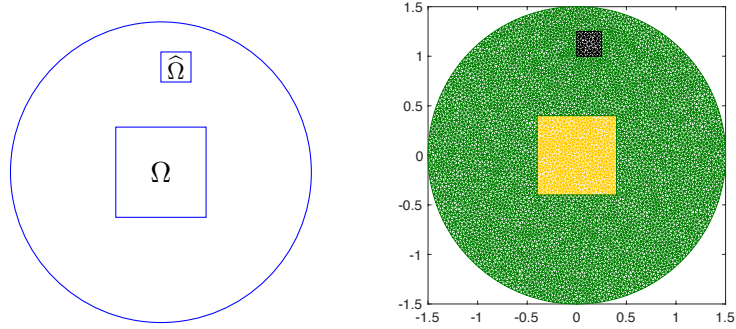
\begin{figure}[!t]
		\centering
			\centering
			\begin{tikzpicture}
				\draw[blue] (0,0) node[black] {$\Omega$} circle (2cm);
				\draw[blue] (0,1.2) rectangle (0.4,1.6);
				\draw[blue] (0.6,0.6) rectangle (-0.6,-0.6);
				\node at (0.2,1.4) (nodeS) {\footnotesize $\widehat{\Om}$};
			\end{tikzpicture}%
			\quad\quad\quad
			\includegraphics[width=0.27\textwidth]{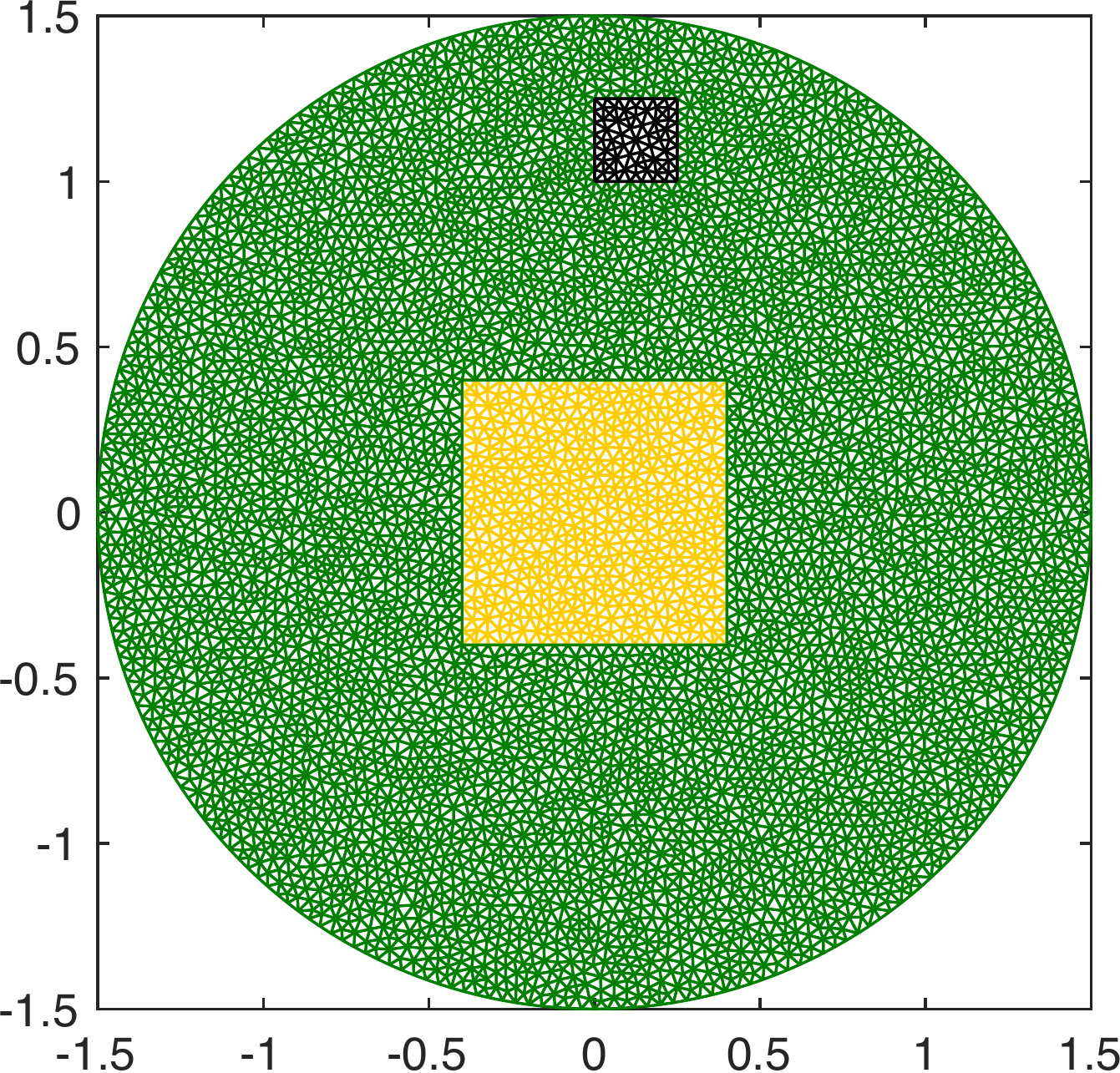}  
			\caption{\label{f:ext_ctrl}
				     \emph{Left:} External optimal control setup for a diffusion process in $\Om$ with control supported in 
				     $\widehat{\Om}$, disjoint from $\Om$. This is different than the classical control
				     approaches where the control must be either inside $\O$ or on the boundary of $\O$.
				     \emph{Right:} A finite element mesh.}			
	\end{figure}

This section provides details on the novel exterior optimal control problem introduced in 
\cite{HAntil_RKhatri_MWarma_2019a,HAntil_DVerma_MWarma_2020a}. Under classical setting,
the control either is on the boundary or in the interior, but this new framework allows 
the control placement away from the domain. This is depicted in Figure~\ref{f:ext_ctrl}. 
Let $Z_D := L^2(\R^N\setminus\O)$, $U_D := L^2(\O)$ and $\lambda \ge 0$ be a constant penalty parameter. 
Then the fundamental optimal control problem amounts to: 
 \begin{subequations}\label{eq:dcp}
 \begin{equation}\label{eq:Jd}
    \min_{(u,z)\in (U_D,Z_D)} J(u) + \frac{\lambda}{2} \|z\|^2_{Z_D} ,
 \end{equation}
 subject to the fractional Dirichlet exterior value problem: Find $u \in U_D$ solving 
 \begin{equation}\label{eq:Sd}
 \begin{cases}
    (-\Delta)^s u &= 0 \quad \mbox{in } \Om, \\
                u &= z \quad \mbox{in } \RR^N\setminus \Om  ,
 \end{cases}                
 \end{equation} 
 and the control constraints 
 \begin{equation}\label{eq:Zd}
    z \in Z_{ad,D} ,
 \end{equation}
 \end{subequations}
 with $Z_{ad,D} \subset Z_D$ being a closed and convex subset.  Well-posedness of this problem follows
 by the standard direct method of calculus of variations \cite[Theorem~4.1]{HAntil_RKhatri_MWarma_2019a}.
 \begin{theorem}\label{thm:exit_dcp}
  Let $Z_{ad,D}$ be a closed and convex subset of $Z_D$. 
  Let either $\lambda > 0$ or $Z_{ad,D}$ bounded and 
  let $J : U_D \rightarrow \mathbb{R}$ be weakly lower-semicontinuous. Then  
  there exists a solution $\bar{z}$ to \eqref{eq:dcp}. If either $J$ is convex and 
  $\lambda > 0$ or $J$ is strictly   convex and $\lambda \ge 0$, then $\bar{z}$ is unique. 
 \end{theorem}
 Moreover, the following first order optimality conditions hold \cite[Theorem~4.3]{HAntil_RKhatri_MWarma_2019a}.
  \begin{theorem}\label{thm:focD}
  Let the assumptions of Theorem~\ref{thm:exit_dcp} hold. Let $\mathcal{Z}$ be an open 
  set in $Z_D$ such that $Z_{ad,D} \subset \mathcal{Z}$. Let 
  $u \mapsto J(u) : U_D \rightarrow \RR$ be continuously Fr\'echet differentiable with 
  $J'(u) \in U_D$. If $\bar{z}$ is a minimizer of \eqref{eq:dcp}
  over $Z_{ad,D}$, then the first order necessary optimality conditions are given by 
  \begin{equation}\label{eq:foc_dcp}
   \left( -\mathcal{N}_s\bar{p}+\lambda\bar{z}, z-\bar{z}\right)_{L^2(\RR^N\setminus\Om)} \ge 0 , \quad 
     \forall  z \in Z_{ad,D},
  \end{equation}
  where $\bar{p} \in \widetilde{H}^{s}(\Om)$ solves the adjoint equation 
  
  \begin{equation}\label{eq:adj_state_dp}
  \begin{cases} 
   (-\Delta)^s \bar{p} = J'(\bar{u}) \quad &\mbox{ in } \Om, \\
               \bar{p} = 0 \quad &\mbox{ in } \RR^N \setminus \Om .
  \end{cases} 
  \end{equation}
  Equivalently we can write \eqref{eq:foc_dcp} as 
  $ \bar{z} = \mathcal{P}_{Z_{ad,D}}\left(\frac{1}{\lambda}\mathcal{N}_s\bar{p}\right)$, 
  where $\mathcal{P}_{Z_{ad,D}}$ is the projection onto 
  the set $Z_{ad,D}$. If $J$ is convex, then \eqref{eq:foc_dcp} is also a sufficient condition. 
 \end{theorem}

As we emphasized earlier, in its current form, \eqref{eq:dcp} requires the approximation of $\mathcal{N}_s$ 
twice. First to compute the very-weak solution (Theorem~\ref{thm:not_soln_ellip}(ii)) and 
second to evaluate the optimality condition \eqref{eq:foc_dcp}. In order to overcome
this issue in the state equation, recall that we had introduced the Robin problem \eqref{eq:Sn2}. 
Also, recall the approximation of the Dirichlet solution by the Robin solution from 
Theorem~\ref{thm:dbc_approx}. In fact, one can approximate the Dirichlet optimal control problem
by the Robin optimal control problem: for $\lambda \ge 0$ the fractional Robin control problem is given by:
 \begin{subequations}\label{eq:ncp_reg}
 \begin{equation}\label{eq:Jn_reg}
    \min_{u\in U_R, z\in Z_R} J(u) + \frac{\lambda}{2} \|z\|^2_{L^2(\Omc)} ,
 \end{equation}
 subject to the regularized exterior value problem (Robin problem): 
 Find $u_n \in U_R$ solving 
 \begin{equation}\label{eq:Sn_reg}
 \begin{cases}
    (-\Delta)^s u &= 0 \quad \mbox{in } \Om \\
    \mathcal{N}_s u  +n \kappa u &= n \kappa z \quad \mbox{in } \RR^N\setminus \Om ,
 \end{cases}                
 \end{equation}
 and the control constraints 
 \begin{equation}\label{eq:Zn_reg}
    z \in Z_{ad,R} .
 \end{equation}
 \end{subequations}
Here $Z_R:=L^2(\Omc)$, $Z_{ad,R}$ is a closed and convex subset of $Z_R$ and 
$U_R := W^{s,2}_{\Om,\kappa}\cap L^2(\Omc)$.
Then as $n \rightarrow \infty$ the Robin problem \eqref{eq:ncp_reg} approximates
the Dirichlet problem \eqref{eq:dcp}, see \cite[Theorem~6.5]{HAntil_RKhatri_MWarma_2019a}
for details. We conclude this section by providing a numerical example, where we solve
\eqref{eq:dcp} by approximating it with \eqref{eq:ncp_reg}. 

We choose our objective function as 
 \[
  j(u,z) = J(u) + \frac{\lambda}{2} \|z\|^2_{L^2(\Omc)} , \quad \mbox{with} \quad 
   J(u) := \frac12 \|u-u_d\|^2_{L^2(\Om)} ,
 \]
and we let $Z_{ad,R} := \{ z\in L^2(\Omc) \;:\; z \ge 0, \ \mbox{a.e. in } \widehat{\Om} \}$ 
where $\widehat{\Om}$ is the support set of the control $z$ that is contained 
in $\widetilde\Om\setminus \Om$. Moreover $u_d : L^2(\Om) \rightarrow \mathbb{R}$ 
is the given data (observations). All the optimization problems below are solved using the projected-BFGS
method with Armijo line search.

Our computational setup is shown in Figure~\ref{f:ext_ctrl}. The centered square region 
is $\Om = (-0.4, 0.4)^2$ and $\widetilde\Om = B_0(3/2)$. The smaller square inside 
$\widetilde\Om \setminus \Om$ is $\widehat\Om$ which is the support of the source/control. 
The right panel in Figure~\ref{f:ext_ctrl} shows a finite element mesh with DoFs = 6103. 

We define $u_d$ as follows. For $z = 1$, we first solve the state equation for 
$\tilde{u}$ \eqref{eq:Sn_reg}. To $\tilde{u},$ we then add a normally 
distributed random noise with mean zero and standard deviation 0.02 to $\tilde{u}$. 
We call the resulting expression as $u_d$. Furthermore, we set 
$\kappa=1$, and $n = 1e5$. 

\begin{figure}[htb]
        \includegraphics[width=0.3\textwidth]{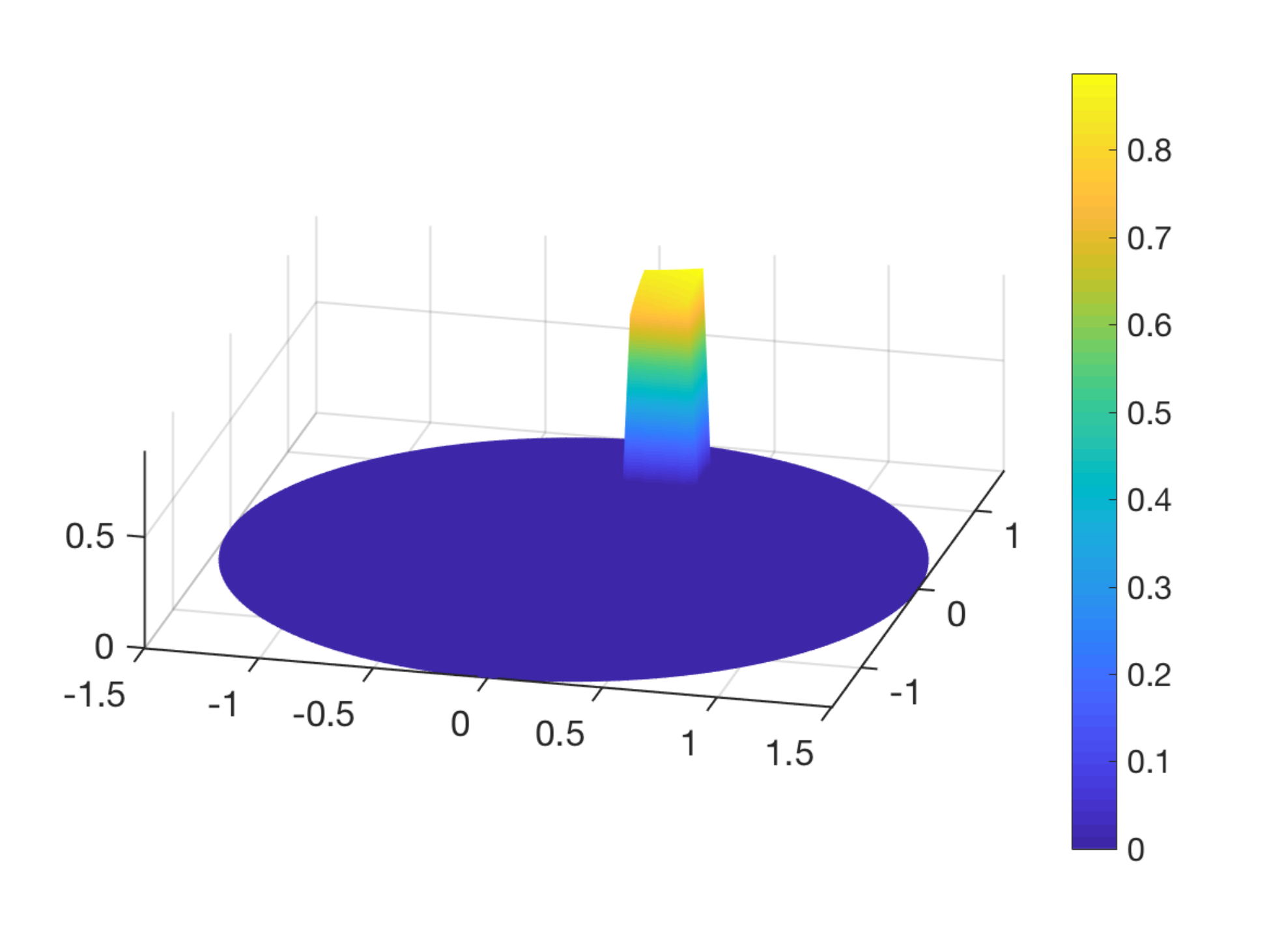}        
        \includegraphics[width=0.3\textwidth]{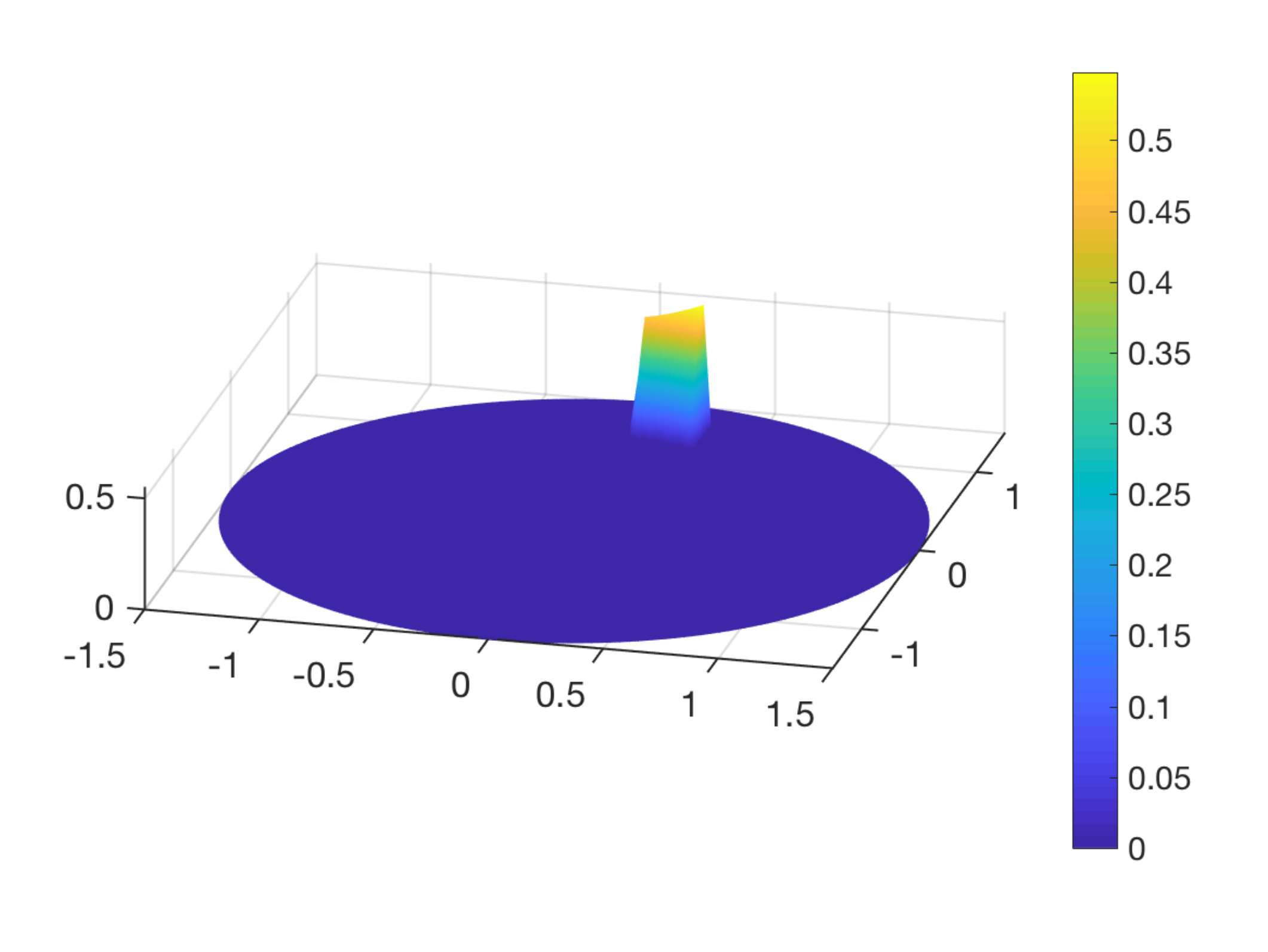}                
        \includegraphics[width=0.3\textwidth]{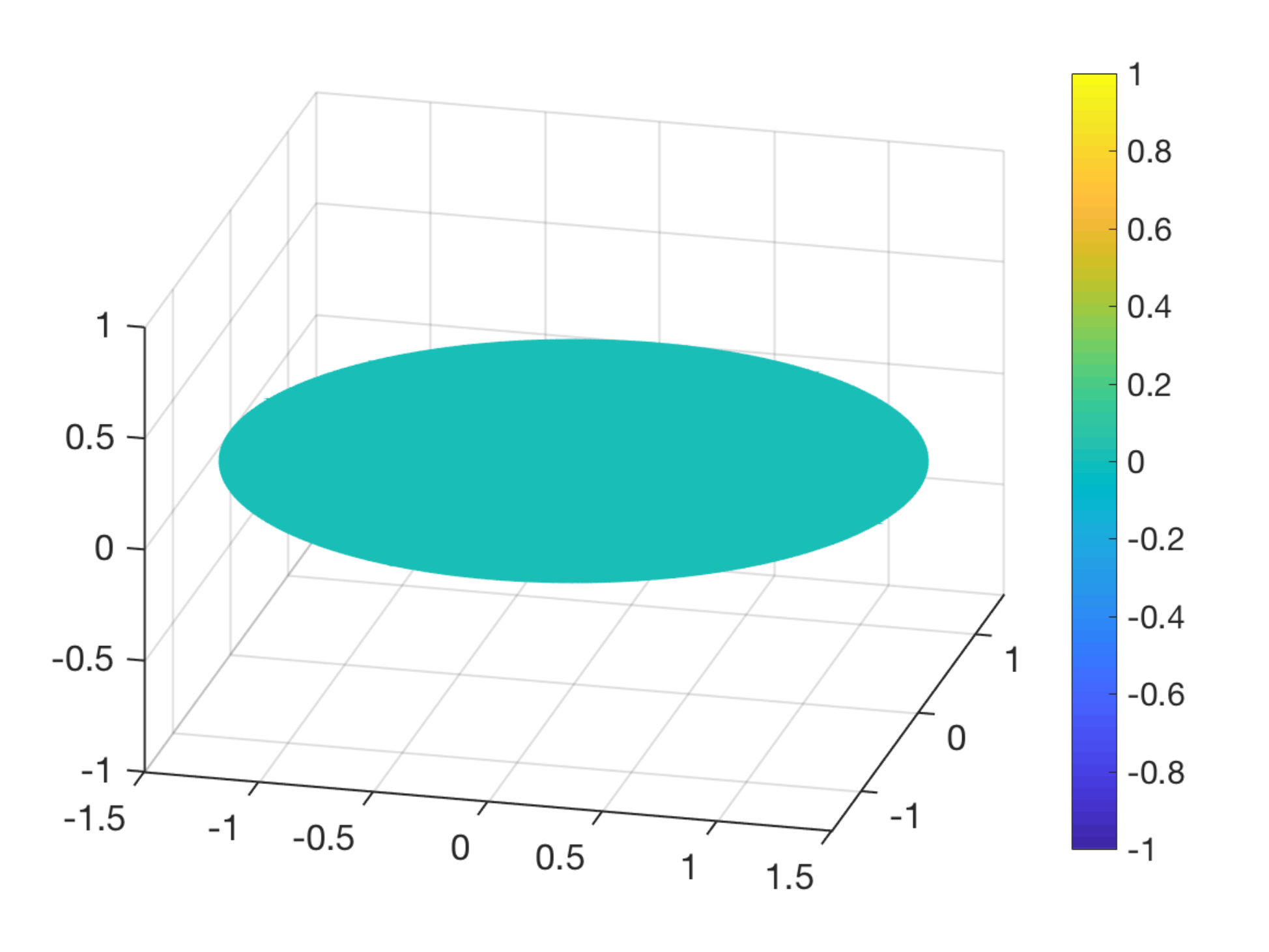}         
        \caption{\label{f:ex4_2}
        The panels show the behavior of $\bar{z}_h$ as we 
        vary the exponent $s$. Left to right: $s = 0.1$, $0.7$, $0.9$. For smaller values of $s$, the recovery
        of $\bar{z}_h$ is quite remarkable. However, for larger values of $s$, $\bar{z}_h 
        \equiv 0$ as expected, the behavior of $\bar{u}_h$ for large $s$ is close to 
        the classical Poisson problem which does not allow external sources/control.}
    \end{figure}

Our goal is then to identify the source/control $\bar{z}_h$. For a fixed $\lambda = 1e-8$, 
Figure~\ref{f:ex4_2} shows the optimal $\bar{z}_h$ for 
$s = 0.1 \ (4)$, $0.7 \ (2)$, $0.9 \ (2)$. The numbers in the
parenthesis denote the total number of iterations that BFGS has taken to achieve a 
stopping tolerance (for the projected gradient) of $1e-7$. Notice that the Armijo line  
search has remained inactive in these cases.
We notice that for large $s$, $\bar{z}_h \equiv 0$. This is 
expected as larger the $s$ is, the closer we are to the classical Poisson problem 
case and we know that we cannot impose the external condition in that case. 

\section{Distributed optimal control of fractional PDEs with state and control constraints} 
\label{s:stateconst}

This section is based on the article \cite{HAntil_TSBrown_DVerma_MWarma_2021a}. 
Given $u_d \in L^2(\O)$ and regularization parameter $\lambda > 0$, consider the state 
constrained optimal control problem
	\begin{equation}\label{eq:st_const} 
	\begin{aligned}
		&\min_{(u,z)\in (U,Z)} \frac12 \| u - u_d\|^2_{L^2(\O)} + \frac{\lambda}{2} \|z\|^2_{L^2(\Om)} \\
		\text{subject to:}& \\
		& (-\Delta)_D^s u = z  \quad \mbox{in } \Om,\qquad u|_{\Om} \in \mathcal{K} \quad \mbox{and} \quad 
		z \in Z_{ad}  .
	\end{aligned}	
	\end{equation}
Recall the definition of $(-\Delta)^s_D$ from \eqref{eq:DLap}. 	
Next, we introduce the relevant function spaces. We let 
    \begin{equation*}
    \begin{aligned}
        Z := L^p(\Om)\;  \mbox{  and }  \;
        U := 
           \Big \{u\in \widetilde{H}^{s}(\Om): ((-\Delta)^s_D)u|_\Om \in L^{p}(\Om)\Big\}.
    \end{aligned}    
    \end{equation*}
Then, $U$ is a Banach space with the graph norm 
$\|u\|_{U}:=\|u\|_{\widetilde{H}^{s}(\Om)} + \|u\|_{C_0(\Om)} 
+ \|(-\Delta)^s_Du\|_{L^p(\Om)}$. Here $p$ is as in \eqref{cond-p},
but in addition we assume that $2 \le p < \infty$. 
We let $Z_{ad} \subset Z$ a nonempty, closed, and convex set and 
$\mathcal{K}$ defined as
    \begin{equation}\label{eq:Ud_1}
       \mathcal{K} := \left\{ w \in C_0(\Omega) \ : \ w(x) \le u_b(x) , 
            \quad \forall x \in \overline{\Om}  \right\}. 
    \end{equation}
Recall that $C_0(\Om)$ is the space of all continuous functions in $\overline\Om$ that vanish
on $\pOm$. Moreover, $u_b \in C(\overline\Om)$ such that $u_b \ge 0$ on $\pOm$. Existence
of solution to \eqref{eq:st_const} then follows from the direct method of calculus of variations 
\cite[Theorem~5.1]{HAntil_DVerma_MWarma_2020b}. 

In addition, under the Slater's constraint qualification
 \cite[Assumption~5.2]{HAntil_DVerma_MWarma_2020b} we can derive the first
order optimality conditions: 
Let $(\bar{u},\bar{z})$ be a solution to the optimization
    problem \eqref{eq:dcp}. Then, there exist a Lagrange multiplier 
    $\bar{\mu} \in (C_0(\Om))^\star = \mathcal{M}(\Om)$ and an adjoint variable $\bar{\vartheta} \in L^{p'}(\Om)$ such that 
    \begin{subequations}\label{eq:st_const_opt}
        \begin{align}
            (-\Delta)^s_D \bar{u} &= \bar{z}  \quad \mbox{in } \Om , \label{eq:a}  \\             
                  \langle \bar{\vartheta},(-\Delta)^s_D v \rangle_{L^{p'}(\Om),L^p(\Om)} 
                &= \left( \bar{u} - u_d ,v \right)_{L^2(\Om)}  
                   + \int_\Om v \; d\bar{\mu}    
                   , && \forall \;v \in U \label{eq:b} \\
                \langle \bar{\vartheta} + J_z(\bar{u},\bar{z}) , 
                    z - \bar{z} \rangle_{L^{p'}(\Om),L^p(\Om)} 
                 &\ge 0 , && \forall \;z \in Z_{ad}  \label{eq:c}  \\
             \bar{\mu} \ge 0, \quad  \bar{u}(x) \le u_b(x) \mbox{ in } \Om, 
             \quad &\mbox{and} \quad \int_{\Om} (u_b - \bar{u})\;d\mu = 0 \label{eq:d} .      
        \end{align}
    \end{subequations}
    
Clearly, it is challenging to try to directly implement \eqref{eq:st_const_opt}. Instead,
we follow the approach from the classical case of $s = 1$ 
\cite{MHinze_RPinnau_MUlbrich_SUlbrich_2009a,MHintermueller_MHinze_2009a}
and consider the so-called Moreau-Yosida regularization. The Moreau-Yosida
regularized optimal control problem is given by 

\begin{subequations}\label{eq:regdcp}
	\begin{equation}\label{eq:regJd}
	\min J^\ga(u,z):=\frac{1}{2}\|u-u_d\|_{L^2{(Q)}}^2 + \frac{\lambda}{2} \|z\|^2_{L^2(Q)} + \frac{1}{2 \ga} \|(\hat{\mu}+\ga (u-u_b) )_+\|^2_{L^2(Q)},
	\end{equation} 
	subject to 
	\begin{equation}
	(-\Delta)^s_D u = z \quad \mbox{in } \O  \quad \mbox{and} \quad \quad z \in Z_{ad}.
	\end{equation} 
\end{subequations}
where $0\leq \hat{\mu}\in L^2(Q)$ is a shift parameter that can be taken to be zero, and $\ga >0$ denotes the regularization parameter. 
Here, $(\cdot)_+$ denotes max$\{0,\cdot\}$. More information about this can be found in \cite{KIto_KKunisch_2008a}. Existence and
uniqueness to \eqref{eq:regdcp} again follows by the standard direct method of calculus of variations. 
Moreover, we have the following first order optimality conditions:
\begin{theorem}\label{thm:optcondRegEll}
	Let $(\bar{u}^\ga,\bar{z}^\ga)$ be a solution to the regularized optimization
	problem \eqref{eq:regdcp}. Then there exists a Lagrange multiplier $\bar\vartheta ^\ga \in \widetilde{H}^s(\O)$ such that 
	\begin{subequations} \label{eq:regEll}
		\begin{align}
		&( -\Delta)^s_D \bar u ^\ga = \bar{z} ^\ga , \quad  &\mbox{in } \Omega,\label{eq:regaEll} \\	
		&(-\Delta)^s_D \bar{\vartheta} ^\ga 
		= \bar{u} ^\ga - u_d + (\hat{\mu} +\ga (\bar u ^\ga - u_b))_+ , \quad &\mbox{in } \Omega ,\label{eq:regbEll} \\
		 & (\bar{\vartheta}^\ga + \lambda \bar{z} ^\ga ,z - \bar{z}^\ga )_{L^2(\Omega)} 
		\ge 0 , \quad &\forall \;z \in Z_{ad} \label{eq:regcEll} .
		\end{align}
	\end{subequations} 
\end{theorem}

In addition, the following approximation result holds.
\begin{proposition}
	Let $(\bar{u},\bar{z})$ solve \eqref{eq:st_const} and $(\bar{u}^\gamma,\bar{z}^\gamma)$ 
	solve \eqref{eq:regdcp}. Then as $\gamma \rightarrow \infty$, we have 
	$\bar{z}^\gamma \rightarrow \bar{z}$ in $L^2(\O)$ and 
	$
		\|(\bar{u}^\gamma - u_b)_+\|_{L^2(\Omega)} = \mathcal O(\gamma^{-1/2}) .
	$ 
	Moreover if $u_b \equiv 0$ in $\O$, then under the additional assumption that $u_d \ge0$,
	we obtain that 
	$
		\|(\bar{u}^\gamma )_+\|_{L^2(\Omega)} = \mathcal O(\gamma^{-1}) 
	$
	as $\gamma\rightarrow \infty$. 
\end{proposition}

For simplicity of presentation, next we consider the case where the control admissible set is
\begin{equation} \label{eq:admis}
	Z_{ad} =  \{ z \in Z : a \leq z \leq b \quad a.e. \; \mbox{in} \; \Omega\}, 
\end{equation}
where $a < b$ are constants. Next, we state the regularity of optimal solutions to \eqref{eq:regdcp}
\begin{proposition}\label{prop:regularity}
Let $\Omega$ be a bounded Lipschitz domain and $(\bar{u}^\gamma, \bar{\vartheta}^\gamma, \bar{z}^\gamma)$ be the solution to \eqref{eq:regEll} for a fixed $\gamma$, then we have 
\[
\bar{u}^\gamma \in H^{\sigma - \varepsilon} (\Omega), \qquad \bar{\vartheta}^\gamma \in H^{\sigma- \varepsilon}(\Omega), \qquad \bar{z}^\gamma \in H^{\tau}(\Omega),
\]
for every $\varepsilon > 0$ where $\sigma = \min\{2s, s+1/2\}$ and $\tau = \min\{1, 2s - \varepsilon\}$. 
\end{proposition}
We discretize the state and adjoint pair using piecewise linear, globally continuous finite elements for a triangulation
$\mathcal{T}_h$ of $\widetilde\Om$ vanishing in the exterior $\widetilde\Om \setminus \Om$. Here 
we have assumed that $\overline\O \subset \widetilde\Om$. We discretize the control using piecewise 
constants, i.e., the control space is 
$
	Z_{h} := \{z_h \in Z : z_h\big|_T \in \mathcal P_0, \ \forall \ T \in \mathcal T_h\},
$
where $\mathcal P_0$ denotes the space of piecewise constants on the triangulation $\mathcal T_h$. 

\begin{figure}[ht]
\centering
\includegraphics[width=0.3\textwidth]{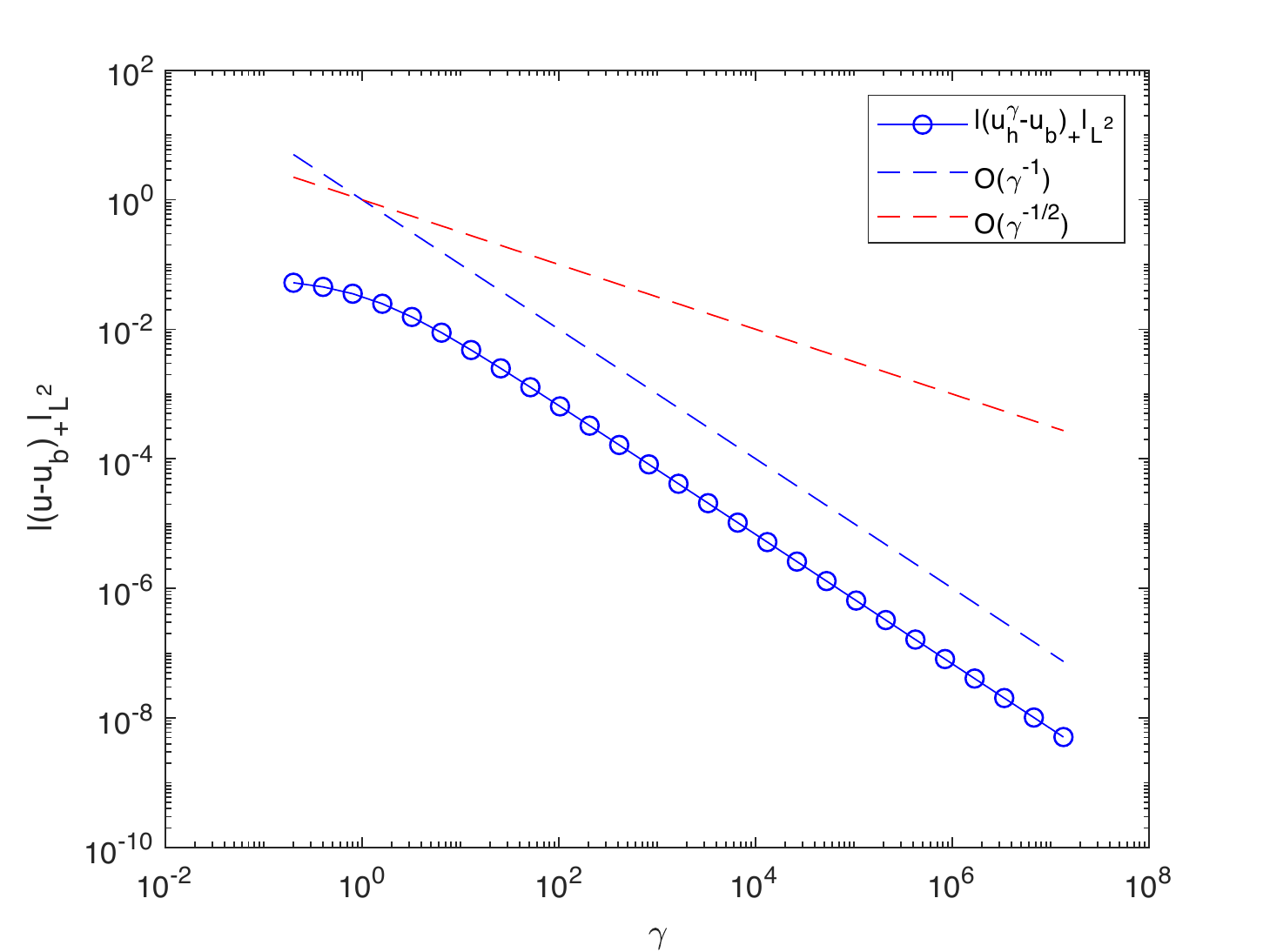} \quad
\includegraphics[width=0.3 \textwidth]{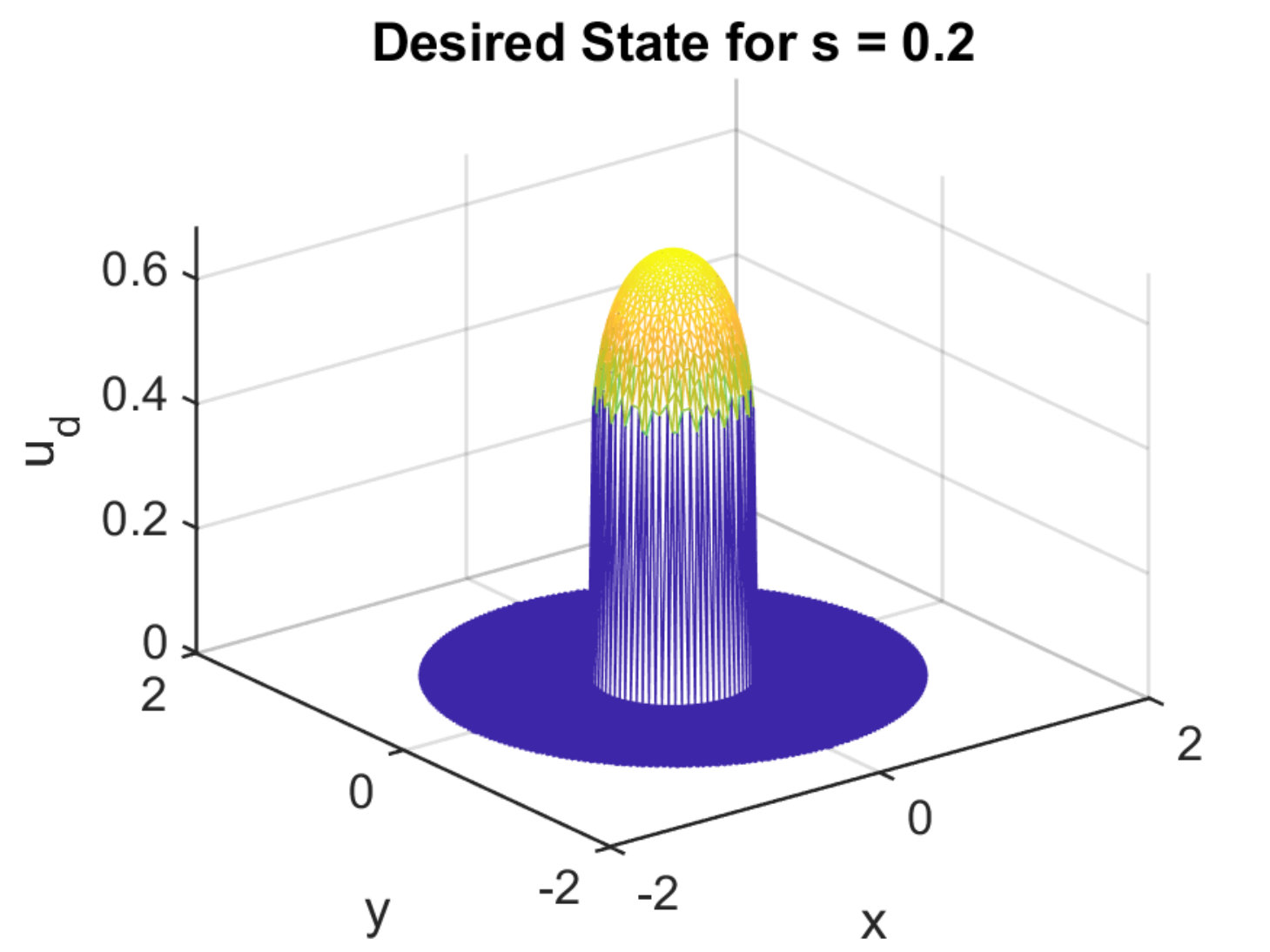} \qquad 
\includegraphics[width=0.3\textwidth]{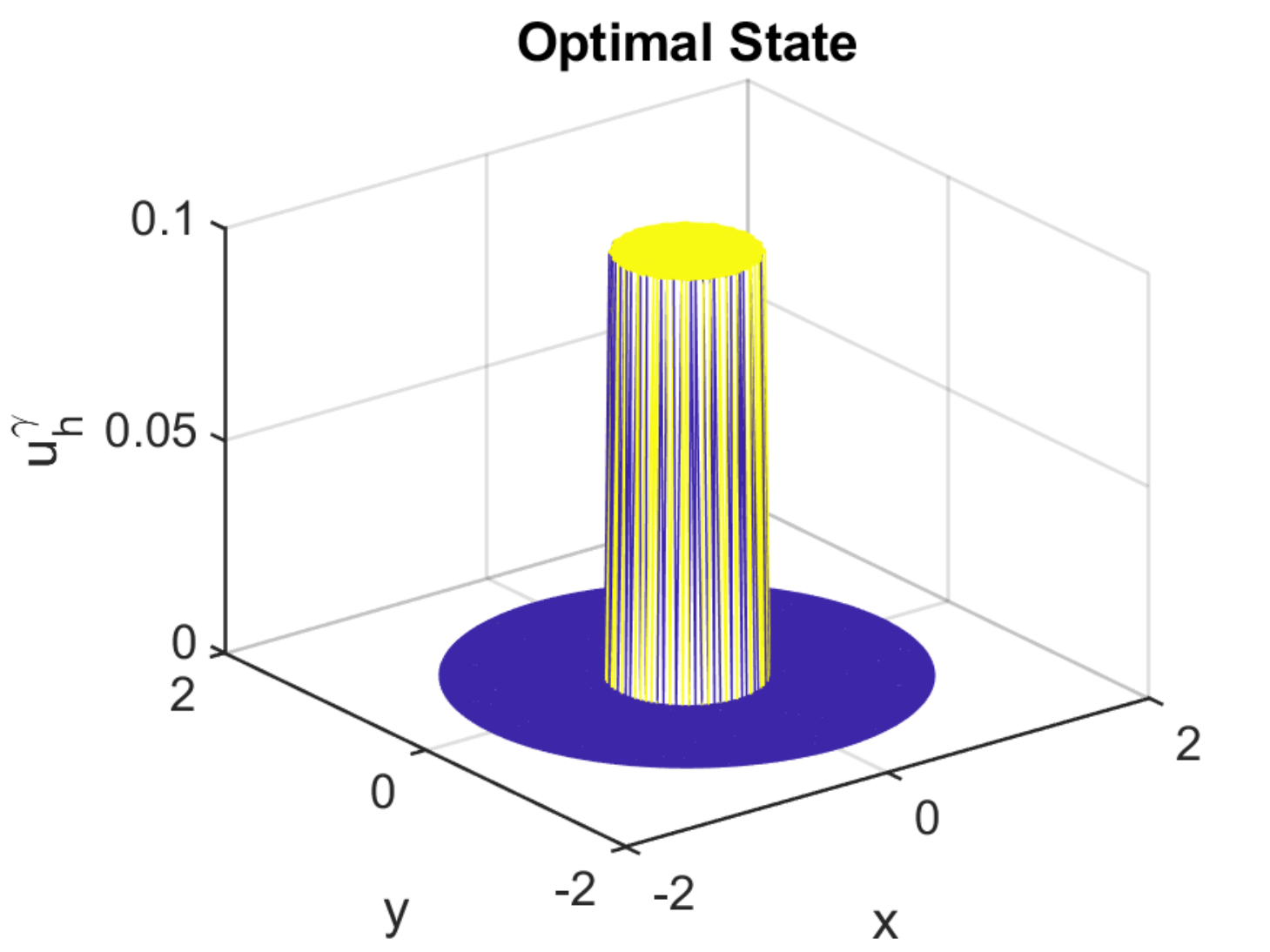}\\
\includegraphics[width=0.3\textwidth]{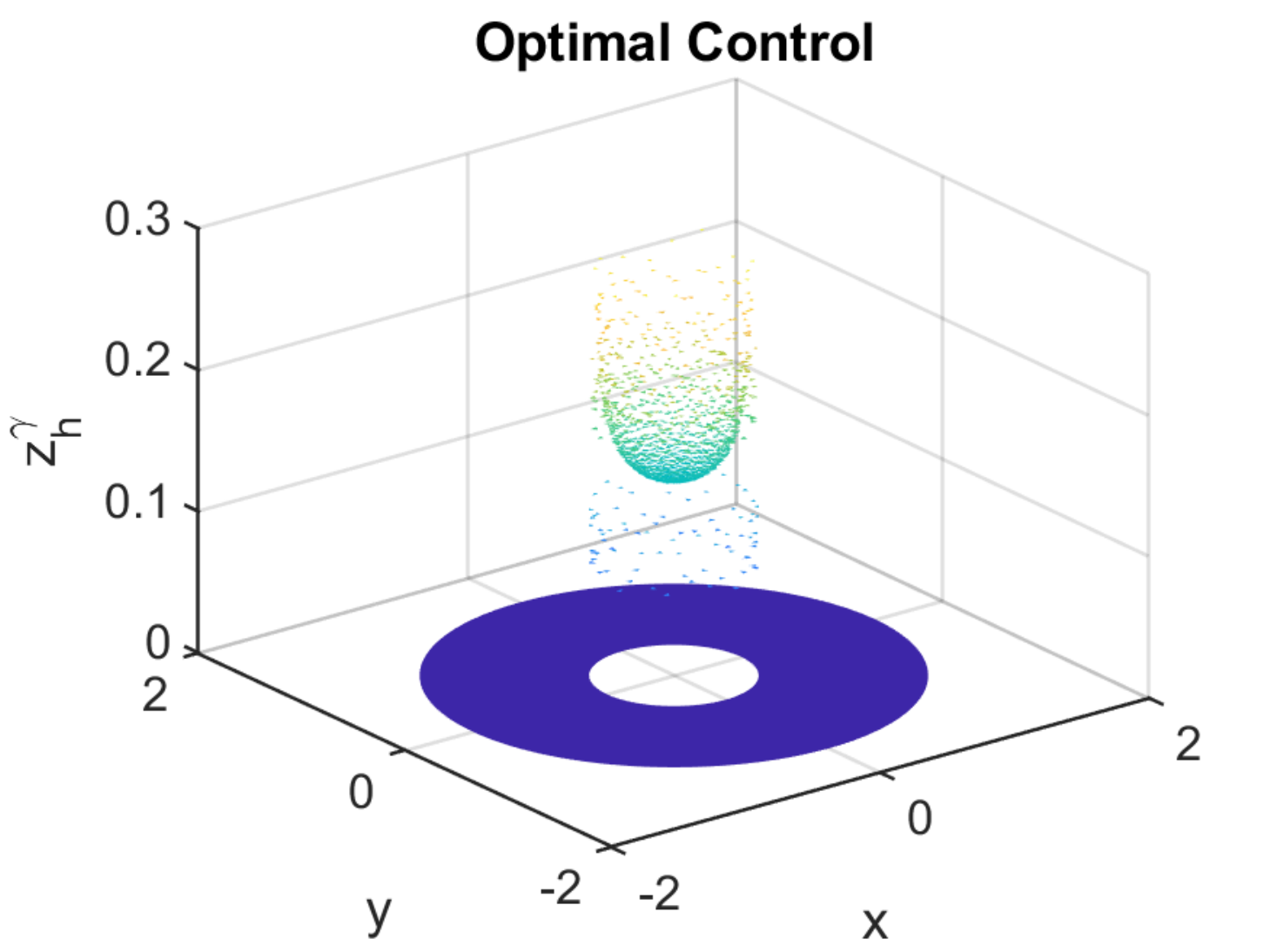} \qquad
\includegraphics[width=0.3\textwidth]{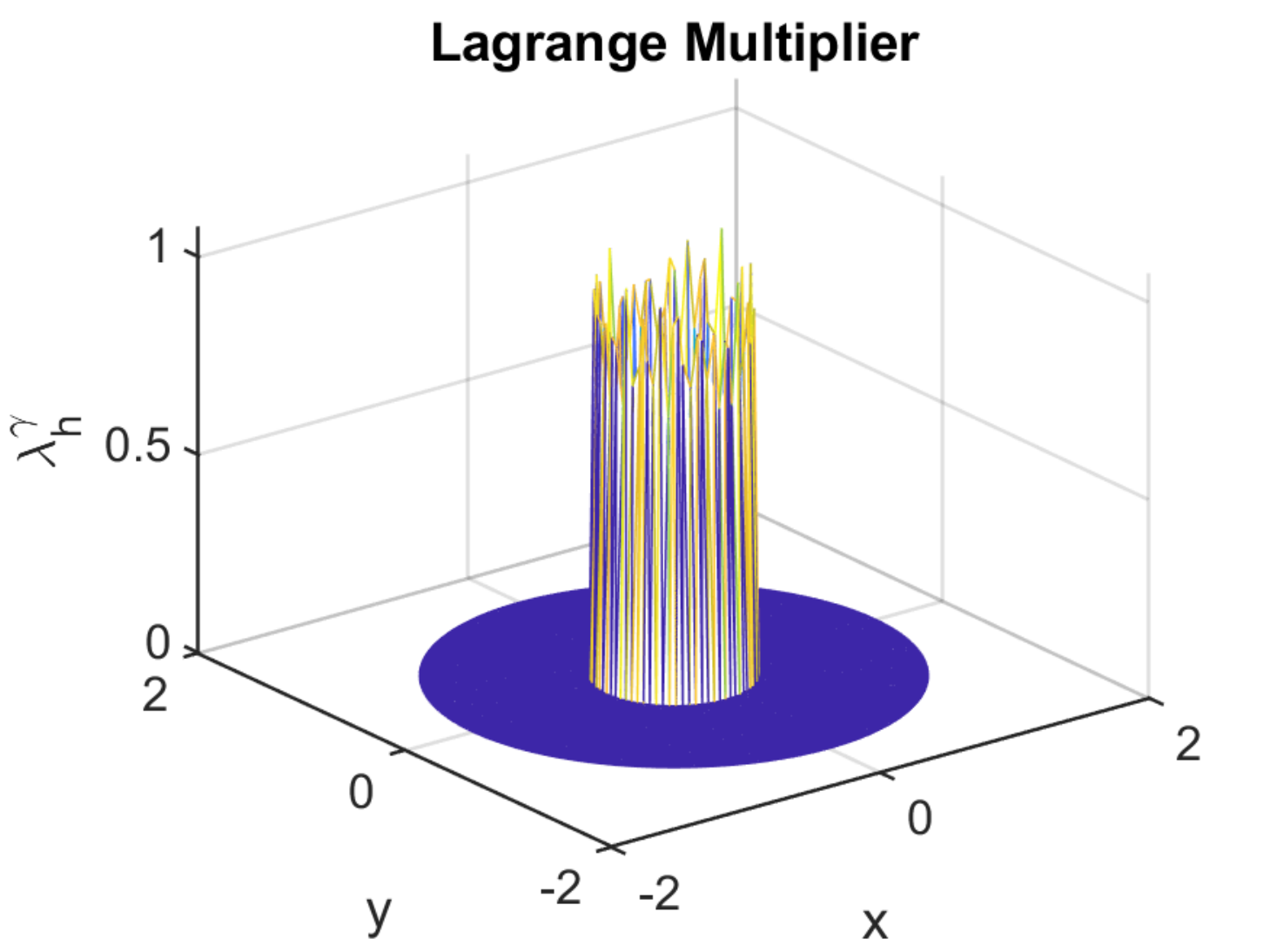}
\caption{Convergence of $\|(u-u_b)_+\|_{L^2(\Omega)}$ as $\gamma$ increases (top left).
 The desired state (top middle), optimal state (top right),  control (bottom left) and Lagrange 
 multiplier (bottom right) for $s = 0.2$ and $\gamma = 419430.4$}
\label{fig:0}
\end{figure}

Then based on the regularity result in Proposition~\ref{prop:regularity}, the following result holds.
\begin{corollary}
\label{cor:err_est}
Let $\bar{z}^\gamma$ and $\bar{z}_h^\gamma$ denote the continuous and discrete optimal controls.
Then the following result holds:
\begin{alignat*}{5}
		\|\bar{z}^\gamma - \bar{z}_h^\gamma\|_{L^2(\Omega)} \leq  \frac{C}{\lambda}\Big(&(h^{2\beta}|\log h|^{2(1+\kappa)} + h^{\beta + s - \varepsilon})  (\|\bar{z}^\gamma\|_{L^2(\Omega)} + \|u_d\|_{L^2(\Omega)} + \|\hat{\mu}\|_{L^2(\Omega)})\\
		& + h^{2\beta}|\log h|^{2(1+\kappa)}(1 + \gamma)\|\bar{z}^\gamma\|_{L^2(\Omega)}
		 	+ h^\tau(1+\gamma + \lambda)|\bar{z}^\gamma|_{H^{\tau}(\Omega)}\Big),
	\end{alignat*}
for every $\varepsilon>0$, where $\beta = \min\{s, 1/2\}, \tau = \min\{1, 2s- \varepsilon\}$, $\kappa = 1$ if $s = 1/2$ 
and zero otherwise.  
\end{corollary}
%
Next, we provide a numerical example. We let 
	$
		u_d (x,y) = \frac{2^{-2s}}{\Gamma(1+s)^2}(1/4-(x^2 + y^2)_+)^s ,
	$
and we let $u_b = 0.1$. 
In the left panel in Figure \ref{fig:0} we show the convergence of $\|(u-u_b)_+\|_{L^2(\Omega)}$ for $s = 0.4$ 
on a mesh with 24155 number of nodes and 48468 number of elements as $\gamma$ increases.  We observe 
a convergence of $\mathcal O(\gamma^{-1})$ which 
is better than expected. However, this has also been documented in the literature (when $s=1$) and it can be 
rigorously established when  $u_b = 0$ and $u_d \ge 0$, see \cite{MHintermueller_MHinze_2009a}. 

In Figure \ref{fig:0} we also show the optimal state, control, and Lagrange multiplier for $s = 0.2$ and 
$\gamma = 419430.4$. We note that the control is a piecewise constant on the mesh.  The optimal state 
in Figure \ref{fig:0} appears to be cleanly cut off at $u_b = 0.1$, complying with the state constraint and 
resulting in a cylindrical profile. Moreover, notice that the Lagrange multiplier corresponding to the inequality 
constraint is a measure (bottom right panel) as expected.  

\section{Fractional deep neural networks -- FDNNs} 
\label{s:time}

This section is based on the article \cite{HAntil_HCElman_AOnwunta_DVerma_2021a}.
Consider a generic (discrete) parameterized PDE 
\begin{equation}
	\label{disc1}
	F(\bu(\bxi);\bxi) = 0, 
\end{equation}
where $\bxi \in \mathcal{P}\subset\mathbb{R}^{N_{\bxi}}$  and $\bu(\bxi)\in \mathcal{U}\subset \mathbb{R}^{N_x}$   
represent  a fixed parameter and the corresponding  solution of the PDE, respectively. Moreover, 
$\mathcal{P}$ denotes the parameter domain and 
 $\mathcal{U}$ the solution manifold. For a fixed parameter $\bxi \in \mathcal{P}$, we seek  
 the solution $\bu(\bxi)\in \mathcal{U}$. In other words, we have the functional relation given 
 by the parameter-to-solution map 
\begin{equation}
	\label{eq:Phi}
 	\bxi\mapsto\Phi(\bxi)\equiv \bu(\bxi).
\end{equation}
Parameter-dependent PDEs of the form  \eqref{disc1} arise in several areas of  computational sciences and engineering. Typical examples include  Navier-Stokes equations  (with Reynolds number as the parameter) \cite{HCElman_DJSilvester_AJWathen_2005a}, and   Boussinesq equations (with Grashof or Prandtl numbers as the parameters)  \cite{AQuarteroni_AManzoni_FNegri_2016a}, etc.

{In real-world applications, solutions of (\ref{disc1}) are required for many parameter values and $N_{\bxi}$ is  often very large; thus, the associated  computational  complexity is enormous.}Besides, a relatively large $N_x$ (due to fine mesh in the discretization of the PDE) yields large (nonlinear) algebraic systems  which are computationally expensive to solve and may also lead to huge storage requirements. This is, for instance, the case in Bayesian inverse problems governed by PDEs where several forward solves are required to adequately sample posterior distributions through  MCMC-type schemes \cite{ABeskos_MGirolami_SLan_PEFarrell_AMStuart_2017a, JMartin_LCWilcox_CBurstedde_OGhattas_2012a}.
Due to the aforementioned challenges, it is a reasonable computational practice to replace the high-fidelity model by a  surrogate model which is relatively easy to evaluate.  

There are two major classes of  surrogate models in the literature: the reduced-order models (ROMs), see e.g., \cite{MBarrault_YMaday_NCNguyen_ATPatera2004a,JSHesthaven_GRozza_BStamm_2016a,AQuarteroni_AManzoni_FNegri_2016a,HAntil_MHeinkenschloss_DCSorensen_2013a}  and the deep neural network  models (DNNs), see e.g.,\cite{YBengio_IGoodfellow_ACourville_2017a, MRaissi_PPerdikaris_GEKarniadakis_2019a}.  A key feature of the ROMs is that they use the so-called {\it offline-online paradigm}. 
The main thrust of the offline step is the construction of a  low-dimensional approximation to the solution space; this approximation is generally known as the {\it reduced basis}. Depending on the problem under consideration, the offline step can be compuationally demanding, although the expense incurred is a one time cost. In the  online step, one then uses the reduced basis to solve a smaller reduced problem. The resulting reduced solution  accurately approximates the solution of the original problem. Typical examples of ROMs include the reduced basis method \cite{AQuarteroni_AManzoni_FNegri_2016a}, proper orthogonal decomposition \cite{JSHesthaven_GRozza_BStamm_2016a, AQuarteroni_AManzoni_FNegri_2016a} and the discrete empirical
interpolation method  (DEIM) and its variants \cite{MBarrault_YMaday_NCNguyen_ATPatera2004a,SChaturantabut_DCSorensen_2010a,HAntil_MHeinkenschloss_DCSorensen_2013a}.

Deep neural network (DNN) models constitute another class of surrogate models which are well-known for their high approximation capabilities. The basic idea of DNNs is to approximate  multivariate functions through a set of layers of increasing complexity \cite{YBengio_IGoodfellow_ACourville_2017a}. As in the case of ROMs, we also note here that using DNNs also involves some offline cost, which is incurred in training the network. Examples of DNNs for surrogate modeling include Residual Neural Networks (e.g. ResNet) \cite{KHe_XZhang_SRen_JSun_2016a,EHaber_LRuthotto_2018a}, physics-informed neural network (PINNs) \cite{MRaissi_PPerdikaris_GEKarniadakis_2019a} and fractional DNN \cite{HAntil_RKhatri_RLohner_DVerma_2020a}.
 
Despite the fact that, in the online phase of ROMs, one essentially has to do reduced system solves, note   that these systems can still be ill-conditioned and  highly intrusive especially for nonlinear problems. On the other hand, the DNN approach can be fully non-intrusive, which is essential for legacy codes. Undoubtedly, rigorous error estimates for ROMs under various assumptions have been well studied \cite{AQuarteroni_AManzoni_FNegri_2016a}; however,
we like the advantage of DNN  being nonintrusive, but recognize that error analysis is not yet as strong \cite{IDaubechies_RDeVore_SFoucart_BHanin_GPetrova_2019a}.

Next, we provide a description of the fractional derivative based DNN approach introduced in 
\cite{HAntil_HCElman_AOnwunta_DVerma_2021a} and  apply it to a Bayesian inverse  
problem. The idea of a DNN is to approximate the input-to-output map $\Phi:\mathbb{R}^{N_{\bxi}}  \rightarrow \mathbb{R}^{N_x}$
by a surrogate $\widehat{\Phi}$ which is the output of a DNN. 

For sufficiently large   $N_s\in \mathbb{N}$,  suppose that   $\mathbb{E}:=\{\bxi_1,\bxi_2,\cdots,\bxi_{N_s}\}$ is a set of parameter samples with $\bxi_i\in\mathbb{R}^{N_{\bxi}}$, and  
$\mathbb{S} := \{\bu(\bxi_1),\bu(\bxi_2),\cdots,\bu(\bxi_{N_s})\}$ the corresponding snapshots (solutions of the model (\ref{disc1}), with $\bu_i := \bu(\bxi_i) \in\mathbb{R}^{N_x}$). Here, we assume that $\mbox{span}\{\bu(\bxi_1),\bu(\bxi_2),\cdots,\bu(\bxi_{N_s})\}$  sufficiently approximates the space of all possible  solutions of (\ref{disc1}).

The idea of a DNN is to use the parameters $\bxi_j$ as an input to the DNN and try to match the output of DNN $\widehat{\Phi}(\bxi_j;\bth)$ 
with the vectors $\bu_j$. Moreover, $\bth =\{W_j,{\bf b}_j\},$ are the unknown parameters in the DNN that need to be learned. This learning
problem can be cast as an optimal control problem
	\begin{equation}
		\label{lossf}
		\min_{\bth = \{W_j, b_j \}} 
		\mathcal{J}(\bth;\bxi, \bu)= \frac{1}{2N_s}\sum_{j=1}^{N_s}\| \hat{\Phi}(\bxi_j;\bth)-\bu_j \|^2_2 + \frac{\lambda}{2}||\bth||_2^2,
	\end{equation}
subject to the fractional Deep Neural Network as constraint
\begin{align}\label{eq:fRNN} 
\phi_1 &= \sigma(W_0 \phi_0 + {\bf b}_0); \quad \phi_0=\bxi, \nonumber \\
\phi_{j} &= \phi_{j-1} - \sum_{k=0}^{j-2}   a_{j-1-k}\: \left(\phi_{k+1} - \phi_k \right)  
		+ \tau^{\gamma}\Gamma(2-\gamma) \ \sigma(W_{j-1}\phi_{j-1} + {\bf b}_{j-1}),\\& &  2\leq j\leq L-1 ,\nonumber \\
(\widehat{\Phi} := ) \ \phi_L   &= W_{L-1} \phi_{L-1} , \nonumber
\end{align}
where, recall that, the middle equation in \eqref{eq:fRNN} is the $L^1$-time discretization of a fractional differential 
equation (cf.~\eqref{Lcaputo} and \eqref{eq:disc_cCaputo}). As pointed out in \cite{HAntil_RKhatri_RLohner_DVerma_2020a},
designing the DNN solution algorithms at the continuous level has the appealing advantage of architecture independence; 
in other words, the number of optimization iterations remains the same even if the number of layers is increased. 

Also, as noted in \cite{HAntil_HCElman_AOnwunta_DVerma_2021a}
this fractional DNN allows connectivity among all the layers unlike standard DNNs. This passage of historic information of 
input and gradients across all the subsequent layers allows one to overcome the vanishing gradient issue as discussed 
in \cite{HAntil_RKhatri_RLohner_DVerma_2020a} for classification problems. 

The learning problem \eqref{eq:fRNN} can be solved using adjoint based approach as discussed in the previous sections
and by applying gradient based optimization methods, such as BFGS. We refer to 
\cite{HAntil_RKhatri_RLohner_DVerma_2020a,HAntil_HCElman_AOnwunta_DVerma_2021a} for the details. We close
this section by providing a numerical example.

\begin{figure}[htb]
	\centering
	\includegraphics[width=0.3\textwidth]{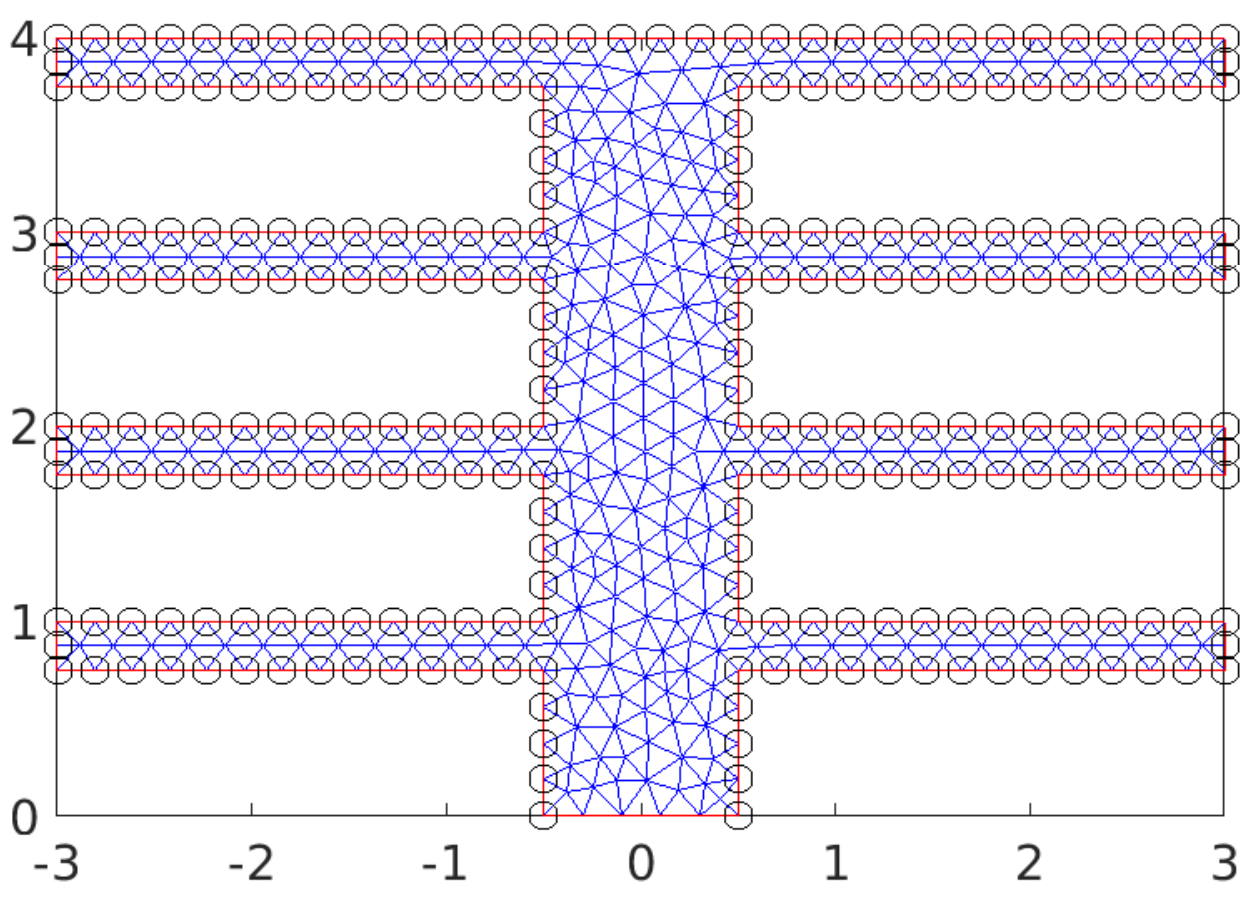}
	\includegraphics[width=0.3\textwidth]{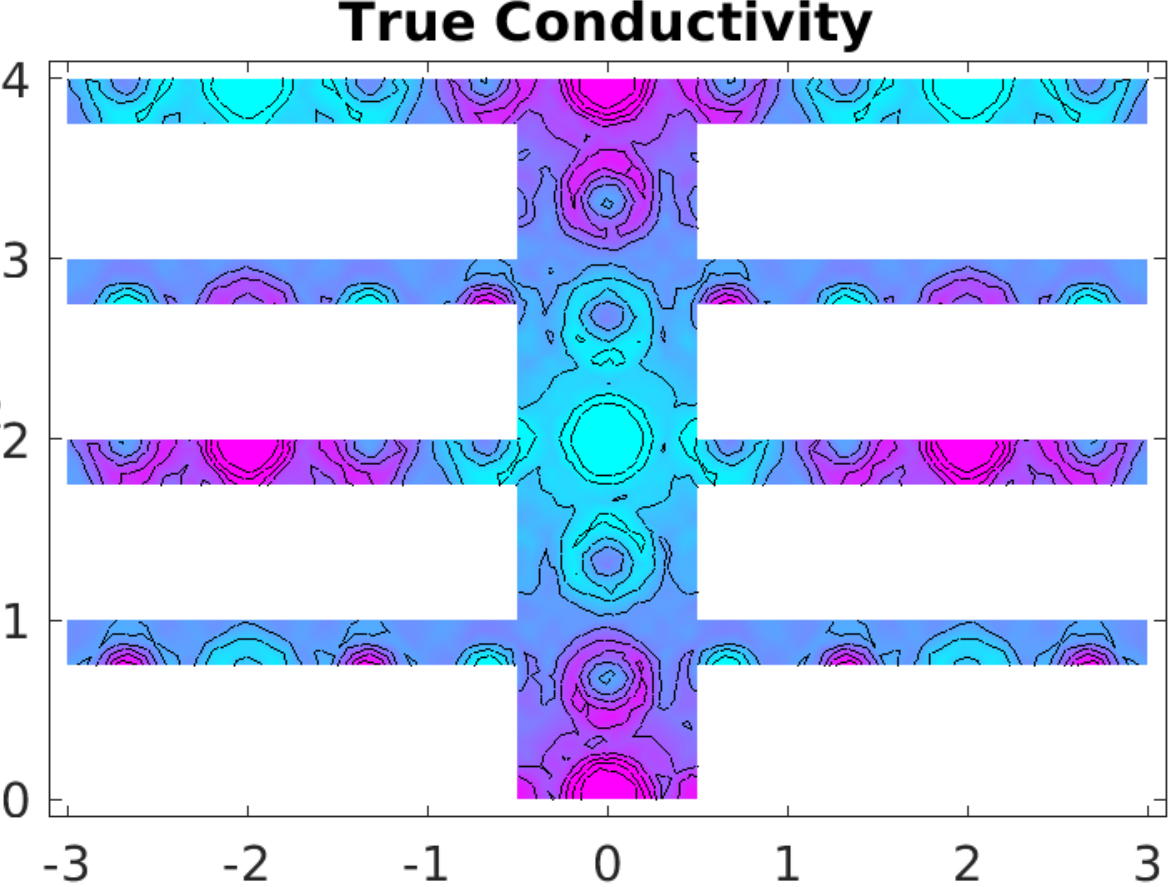}		
	\includegraphics[width=0.3\textwidth]{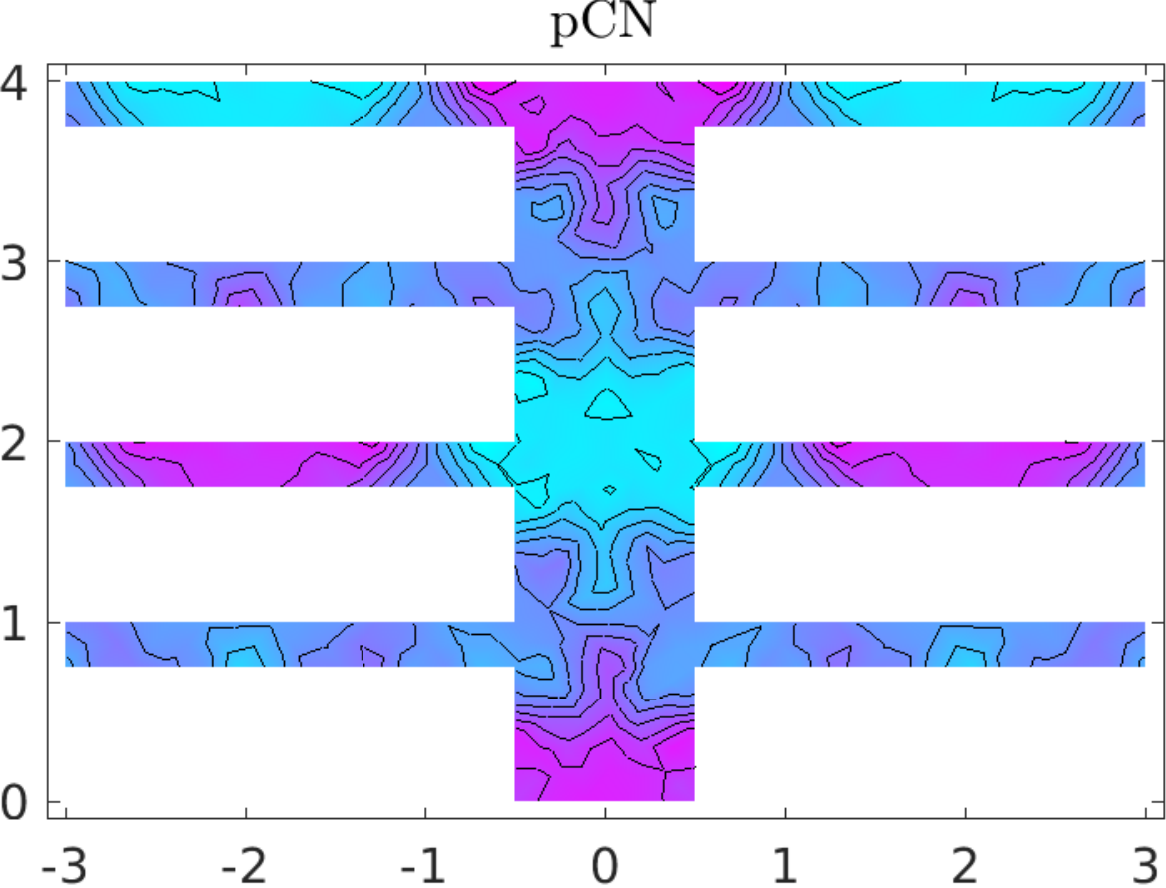}	
	\caption{ \emph{Thermal fin problem:} Given a diffusion equation with coefficient $e^{\bxi(x)}$, 
	the goal is to reconstruct $u$, from noisy measurements on the boundary (marked by circles). 
	\emph{Middle:} True conductivity field. \emph{Right:} The posterior mean estimates 
	obtained by preconditioned Crank-Nicolson (pCN) MCMC method using the fractional DNN as 
	surrogate. The acceptance rate is comparable 
	to \cite{ABeskos_MGirolami_SLan_PEFarrell_AMStuart_2017a}. 
	The figure has been reproduced from \cite{HAntil_HCElman_AOnwunta_DVerma_2021a}.
	}
		\label{f:UQ}
\end{figure}

\medskip 
\noindent
{\bf Thermal fin problem:}
%
Consider the following parameterized system  
\begin{align}
 \label{soc}
 \begin{aligned} 
    -\mbox{div } (e^{\bxi(x)}\nabla u )  &=0,  \quad  &&\mbox{in } \Omega , \\
   (e^{\bxi(x)}\nabla u )\cdot \boldsymbol{\nu} + 0.1 u &= 0, \quad &&\mbox{on }  \partial\Omega \setminus\Gamma,\,\\
   (e^{\bxi(x)}\nabla u )\cdot {\bf n} &= 1, \quad &&\mbox{on } \Gamma=(-0.5,  0.5) \times \{ 0 \},
\end{aligned}
\end{align}
which represents a forward model for  heat conduction over the non-convex domain $\Omega$ as shown in 
Figure~\ref{f:UQ} (left). When the heat conductivity function $e^{\bxi({\bf x})}$ is known, the forward problem 
can be solved for $u$. The goal of this example is to infer $100$ unknown parameters
$\bxi$ from  $262$  noisy observations of $u$. Figure~\ref{f:UQ} (left) also shows the location of the  observations on 
the boundary $\partial\Omega \setminus\Gamma$. 

The fDNN is trained by first computing $N_s = 900$ solution snapshots corresponding to 
900 parameters $\bxi \in \mathbb{R}^{100}$, drawn using Latin hypercube sampling. Next we compute the
SVD of $\mathbb{S}$ and train the network on this reduced space. Figure~\ref{f:UQ} shows the true 
conductivity (middle) and the posterior mean estimates obtained using standard MCMC approach. 

Table~\ref{acc_time} below shows the average acceptance rates for these models and the 
computational times required to perform the MCMC simulation for different variants of the 
MCMC algorithm
(preconditioned Cranck-Nicolson method (pCN), the infinite variants of Riemannian manifold 
Metropolis-adjusted Langevin algorithm ($\infty$-MALA)  and  Hamiltonian Monte-Carlo 
($\infty$-HMC) algorithms as presented in \cite{ABeskos_MGirolami_SLan_PEFarrell_AMStuart_2017a}),
using both fDNN surrogate computations and full-order discrete PDE solution.
The acceptance rate is clearly comparable for the fDNN surrogate model and the full order forward
discrete PDE solver. In addition, a reduction of 90\% in computational time is observed for the
fDNN approach. The costs of the offline computations for fDNN were 
$63.2$ seconds to generate the data used for the fractional DNN models and  $33.5$ seconds to 
train the network (with $1600$ BFGS iterations); a total of $99.7$ seconds.

 \begin{table}[h!]
 \centering
\begin{tabular}{l|lllll}
\hline
Model    &  pCN        &  $\infty$-MALA          &  $\infty$-HMC &  \\
\hline
\hline
Acc. Rate (fDNN)  &  $0.67$  &   $0.67$    & $0.79$  & \\
\hline
Acc. Rate (Full)  &  $0.66$  &   $0.70$    & $0.75$  & \\
\hline
\hline
CPU time  (fDNN) &  $16.46$  &    $ 98.0$    & $228.4$   & \\
\hline
CPU time  (Full) &  $157.8$  &    $ 958.9$    & $2585.3$   & \\
\hline
\end{tabular}
\caption{Acceptance rates and computational times  
needed 
 to solve the inverse problem by pCN, $\infty$-MALA and $\infty$-HMC algorithms together with fDNN  
 and  full forward  models.
}
\label{acc_time}
\end{table}


\section{Some open problems}
\label{s:open}

There are still many open questions in this paradigm. We list a few of them below:
\begin{enumerate}
	\item Theorem 3.2 in part (iv) assumed the domain $\Om$ to fulfill exterior cone condition. It remains
	        open if such a result holds true in case $\Om$ is only Lipschitz continuous. Also, note that (iii)-(iv) shows 
		boundedness and continuity of solution when the exterior datum $g \equiv 0$. However, it 
		remains open to prove these results when $g \neq 0$.  That is,  given a function $g\in C(\mathbb R^N\setminus\Omega)$ and $f\equiv 0$ in $\Omega$,  under what conditions on $\Omega$, there would exist an $s$-harmonic function $u\in C(\mathbb R^N)$ satisfying \eqref{eq:ell_d_gen}. The classical case $s=1$ has been resolved by Wiener several decades ago, but the fractional case is more challenging and remains an open problem.
		
	\item Section~\ref{s:nhbc} discusses a numerical method to approximate the numerical method to 
		approximate the non-homogeneous exterior Dirichlet problem by the Robin problem. But a 
		complete numerical analysis of this method is still open. In addition, a complete numerical 
		analysis of the Robin problem in the full generality considered here is still open as well. 
		
	\item Section~\ref{s:ext} discusses exterior optimal control problems, but numerical analysis of 
		 these problems are fully open. 
		 
	\item Section~\ref{s:stateconst} discusses state constrained problems and provides finite element 
		$\gamma$-dependent error estimates in Corollary~\ref{cor:err_est}. Such estimates which are 
		$\gamma$-independent are still open. Note that these questions are not just limited to 
		elliptic problems, but are equally relevant to parabolic problems as well. 
\end{enumerate}


\bibliographystyle{plain}
\bibliography{references}

\end{document}